\documentclass[10pt]{amsart}

\usepackage{amsfonts, amstext, amsmath, amsthm, amscd, amssymb}
\usepackage{epsfig, graphics, psfrag}
\usepackage{color}

\renewcommand{\setminus}{{\smallsetminus}}

\newcommand{\cross}{{\times}} 
\newcommand{\bdy}{{\partial}} 

\newcommand{\len}{{\mathrm{Length}}}
\newcommand{\genus}{{\mathrm{genus}}}
\newcommand{\half}{{\frac{1}{2}}}
\newcommand{\lgeom}{\ell_g} 
\newcommand{\lcomb}{\ell_c} 
\newcommand{\abs}[1]{{\left\vert #1 \right\vert}}
\renewcommand{\inf}[1]{{\mathrm{inf} \left\{ #1 \right\} }}

\newcommand{\tild}[1]{{\widetilde{#1}}}

\newcommand{\HH}{{\mathbb{H}}}
\newcommand{\RR}{{\mathbb{R}}}
\newcommand{\ZZ}{{\mathbb{Z}}}
\newcommand{\CC}{{\mathbb{C}}}

\newcommand{\s}{{\boldsymbol{s}}}
\newcommand{\w}{{\boldsymbol{w}}}

\theoremstyle{plain}
\newtheorem{theorem}{Theorem}[section]
\newtheorem{corollary}[theorem]{Corollary}
\newtheorem{lemma}[theorem]{Lemma}
\newtheorem{prop}[theorem]{Proposition}

\newtheorem*{no-num-theorem}{Theorem}

\theoremstyle{definition}
\newtheorem{define}[theorem]{Definition}
\newtheorem*{remark}{Remark}
\newtheorem*{notation}{Notation}

\begin{document}

\title{Links with no exceptional surgeries}
\author{David Futer}
\author{Jessica S. Purcell}
\date{\today}

\address{David Futer, Department of Mathematics, Stanford University, 
Stanford, CA 94305}
\email{dfuter@math.stanford.edu}

\address{Jessica S. Purcell, Department of Mathematics, 1 University
  Station C1200, University of Texas at Austin, Austin, TX 78712}
\email{jpurcell@math.utexas.edu}

\begin{abstract}
  
  We show that if a knot admits a prime, twist-reduced diagram with at
  least $4$ twist regions and at least $6$ crossings per twist region,
  then every non-trivial Dehn filling of that knot is hyperbolike.  A
  similar statement holds for links. We prove this using two
  arguments, one geometric and one combinatorial.  The combinatorial
  argument further implies that every link with at least $2$ twist
  regions and at least $6$ crossings per twist region is hyperbolic 
  and gives a lower bound for the genus of a link.
\end{abstract}

\maketitle

\section{Introduction}\label{intro}

Knots and links in $S^3$ are easiest to visualize with a projection
diagram, but computing geometric or topological information directly
from the diagram is often a difficult task. In the very special case
of alternating knots, an alternating diagram reveals a lot of
topological information, including the genus of the knot
\cite{crowell, murasugi-genus} and bounds on the Heegaard
genus of the complement \cite{lack-tunnel}. For alternating knots and
links, one can tell by looking at an alternating diagram whether the
complement is hyperbolic \cite{menasco-alt}, and if it is, compute upper
and lower bounds on the volume \cite{lack-volume}. However, few
results of this sort extend beyond this special class of knots and
links.

In this paper, we prove a mild diagrammatic condition that ensures the
complement of a particular link is hyperbolic, and a slightly stronger
one that ensures all non-trivial Dehn surgeries on the link are {\it
hyperbolike}. 
We also use the combinatorial properties of a diagram to give a 
lower bound on the genus of a link.
To state our results precisely, we will need a few
definitions.

\subsection{Twist regions and reduced diagrams}

A \emph{diagram} $D(K)$ of a knot or link $K \subset S^3$ can be
viewed as a $4$-valent planar graph $G(K)$, with over-under crossing
information at each vertex.

\begin{define}
A \emph{bigon} is a contractible region in the complement of $G(K)$
that has two edges in its boundary. Following Lackenby
\cite{lack-surg, lack-volume},
 define a \emph{twist region} of the knot or link
to be a maximal string of bigons arranged end to end.  A single
crossing adjacent to no bigons is also a twist region.
\end{define}

\begin{figure}[ht]
  \begin{center}
    \includegraphics{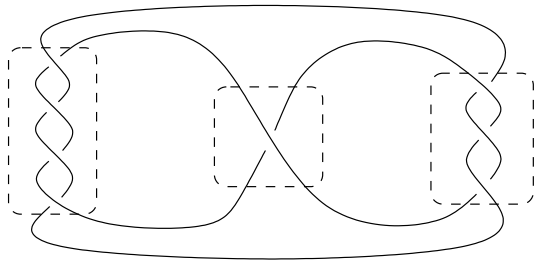}
  \end{center}
  \caption{The above diagram has $3$ twist regions, containing $2$,
    $\frac{1}{2}$, and $1\frac{1}{2}$ twists, respectively.}
  \label{fig:intro_twistnumber}
\end{figure}

We are also concerned with the amount of twisting that occurs in each
twist region.  We will count this either in terms of crossings or in
terms of full twists, where a {\it full twist} of one strand about the
other corresponds to two crossings. See Figure
\ref{fig:intro_twistnumber} for an illustration of these definitions.

\begin{figure}[ht]
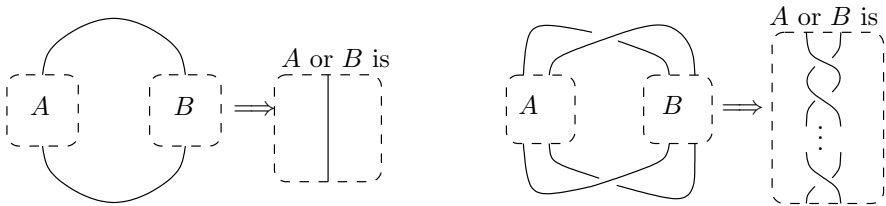

\vspace{-12pt}
  \begin{center}
  \input{Fig_prime_diagr.pstex_t}
  \hspace{.5in}
  \input{Fig_twist_red_diagr.pstex_t}
  \end{center}
  \caption[Prime and twist-reduced diagrams.]{Left: A prime diagram;
  Right: A twist-reduced diagram.}
  \label{fig:prime_twred_diagr}
\end{figure}

\begin{define}
A diagram $D(K)$ of a knot or link $K \subset S^3$ is called
\emph{prime} if for any simple closed curve $\gamma$ in the projection
plane that intersects the graph $G(K)$ transversely in two points in
the interior of edges, $\gamma$ bounds a subdiagram containing no
crossings of the original diagram.  Note this ensures the diagram
contains no monogons.  See Figure \ref{fig:prime_twred_diagr}.
\end{define}
  
Following Lackenby \cite{lack-volume}, we also require the diagram to
be \emph{twist-reduced}.

\begin{define}
A diagram $D(K)$ of a knot or link $K$ is \emph{twist-reduced} if
whenever a simple closed curve $\gamma$ in the projection plane
intersects the graph $G(K)$ transversely in four points in the
interior of edges, with two points adjacent to one crossing and the
other two points adjacent to another crossing, then $\gamma$ bounds a
subdiagram consisting of a (possibly empty) collection of bigons
arranged in a row between these two crossings.  See Figure
\ref{fig:prime_twred_diagr}.
\end{define}

Note that any diagram of a prime knot or link $K$ can be simplified
into a prime, twist-reduced diagram: for if $D(K)$ is a diagram that
fails to be prime, then all crossings on one side of a simple closed
curve $\gamma$ are extraneous and can be removed.  Similarly, if
$D(K)$ is not twist-reduced, then a series of flypes will amalgamate
the two twist regions adjacent to a curve $\gamma$ into a single
region, reducing the number of twist regions.

These definitions allow us to state our first result.

\begin{theorem}\label{hyp-link}
Let $K \subset S^3$ be a link with a prime, twist-reduced diagram
$D(K)$. Assume that $D(K)$ has at least two twist regions (i.e. that
$K$ is not a closed $2$-braid). If every twist region of $D(K)$
contains at least $6$ crossings, then $K$ is hyperbolic.
\end{theorem}

This result could be viewed as an extension of Menasco's theorem
\cite{menasco-alt}, which holds that a prime, non-split alternating
link is hyperbolic whenever it is not a closed $2$-braid. For
alternating links, Menasco doesn't need any assumption on the number
of crossings per twist region; to rule out non-hyperbolic links in
general, some such assumption is necessary.

\smallskip

Recall that the \emph{genus} of a link $K \subset S^3$ is the smallest genus of an incompressible, orientable surface $S \subset S^3$ whose boundary is $K$. We can use prime, twist-reduced diagrams to give a lower bound on the genus.

\begin{theorem}\label{genus-bound}
Let $K \subset S^3$ be a link of $k$ components with a prime, twist-reduced diagram $D(K)$. If $D(K)$ has $t \geq 2$ twist regions and at least $6$ crossings in each twist region, then
$$\genus(K) \geq \left\lceil 1 + \frac{t}{6} - \frac{k}{2} \right\rceil ,$$
where $\lceil \cdot \rceil$ is the ceiling function that rounds up to the nearest integer.
\end{theorem}

Crowell \cite{crowell} and Murasugi \cite{murasugi-genus} have independently proved that the genus of an alternating link is equal to half the degree of its Alexander polynomial, and Gabai \cite{gabai-arborescent} gave an algorithm to compute the genus of an arborescent link. The advantage of Theorem \ref{genus-bound} is that it works for general links and, in fact, gives the exact value for certain families of links.

\subsection{Dehn surgery}

Let $M$ be a 3-manifold with torus boundary $\partial M$, and $s$ a
\emph{slope} on $\partial M$, that is, $s$ is an isotopy class of
simple closed curves on $\partial M$. The manifold obtained by gluing
a solid torus $S^1 \times D^2$ to $\partial M$ in such a way that the
slope $s$ bounds a disk in the resulting manifold is called a
\emph{Dehn filling along the slope $s$}, or a \emph{Dehn surgery along
$s$}.  More generally, if $M$ is a 3-manifold with multiple torus
boundary components and along each component we have a slope $s_i$, we
obtain a closed manifold by Dehn filling along these slopes.

Using a basis $\langle \mu, \lambda \rangle$ for the fundamental group
of the torus, slopes on cusps are parameterized by $\mathbb Q \cup
\{\infty\}$.  Thus a slope corresponds to $a/b$ if and only if the
slope is equivalent to $a\mu+b\lambda$.  If $K$ is a knot in $S^3$,
and $M$ is taken to be the {\it exterior} $E(K)$ of a tubular
neighborhood of $K$, then we let $\mu$ correspond to a meridian, and
$\lambda$ to a longitude.  In this case, Dehn filling along a meridian
of $K$, i.e. $1/0$ filling, will always give $S^3$.  This Dehn filling
is called the \emph{trivial filling}.  All other Dehn fillings are
\emph{non-trivial}.

Thurston \cite{thur-notes} has shown that given a hyperbolic manifold
$M$ with cusps, all but finitely many choices of surgery slope on each
component of $\partial M$ yield a closed hyperbolic manifold.  More
recently, Hodgson and Kerckhoff \cite{univ-bounds} showed that if the
surgery slope on each component of $\bdy M$ is longer than a given
universal constant, then the resulting Dehn filled manifold is
hyperbolic.  Using these results, Purcell \cite{purcell:thesis} was
able to show that for sufficiently complicated knots, every nontrivial
Dehn filling is hyperbolic.  However, the required knots are so
complicated that they are difficult to use in practice.

If we weaken the assumption that the resulting manifold be hyperbolic,
we can obtain similar surgery results for much less complicated knots.

\begin{define}\label{hyplike-def}
A closed, orientable $3$-manifold $M$ is \emph{hyperbolike} if
\begin{enumerate}
  \item $M$ is irreducible and atoroidal,
  \item $M$ is not Seifert fibered, and
  \item $\pi_1(M)$ is infinite and word-hyperbolic.
\end{enumerate}
\end{define}

All hyperbolic manifolds are hyperbolike, and Thurston's
Geometrization Conjecture \cite{thur-survey} would imply the converse.

\smallskip

\begin{theorem}\label{main}
Let $K$ be a link in $S^3$ with a prime, twist-reduced diagram
$D(K)$. Suppose that every twist region of $D(K)$ contains at least
$6$ crossings and each component of $K$ passes through at least $7$
twist regions (counted with multiplicity). Then every non-trivial Dehn
filling of all the components of $K$ is hyperbolike.
\end{theorem}

Notice that we obtain this information about Dehn fillings from the
diagram of the link alone.  We need no additional information.  In
fact, our methods also prove that a non-trivial Dehn filling of only
some components of $K$ yields a hyperbolic manifold with boundary.

\begin{corollary}\label{main-knot-cor}
Let $K$ be a knot in $S^3$ with a prime, twist-reduced diagram
$D(K)$. If $D(K)$ has at least $4$ twist regions, and each twist
region contains at least $6$ crossings, then any non-trivial Dehn
filling of $K$ is hyperbolike.
\end{corollary}

The corollary follows from Theorem \ref{main} because if $K$ is a
knot, every twist region contains two strands of $K$. Thus in a
diagram with $4$ twist regions, $K$ passes through a twist region $8$
times. 

In fact, the hypothesis of $4$ twist regions in Corollary
\ref{main-knot-cor} is a sharp bound.  Wu
\cite{wu:toroidal_montesinos} has shown that every pretzel knot 
with $3$ twist regions and at least $2$ crossings per twist region
admits a non-trivial exceptional surgery.  Thus our results assume the
smallest possible number of twist regions.

As for the requirement that each twist region contain at least $6$
crossings, we know that some such requirement is necessary.
It is known that there exist knots with non-trivial exceptional surgeries
that have arbitrarily large volume, hence an arbitrarily high number of
twist regions.  
These have been discovered by Eudave-Mu{\~n}oz and
Luecke \cite{EM-Luecke:volume}, Eudave-Mu\~noz
\cite{eudave:seifert-surgery}, as well as recently by Baker
\cite{baker}.
Thus a high number of twist regions alone is not enough
to rule out exceptional surgeries.  However, at this time the authors
do not know whether the requirement of six crossings per twist region
is sharp.

Another advantage of Theorem \ref{main} is that it gives information
on Dehn fillings without requiring us to restrict our attention to a
particular class of knots or links.  This should be compared to other
known results.  If we restrict to alternating links, Lackenby
\cite{lack-surg} has shown that all non-trivial Dehn surgeries on
alternating knots with at least 9 twist regions are hyperbolike, as
are surgeries on alternating links in which each component passes
through 17 or more twist region.  Wu \cite{wu-knot} proved that all
non-trivial surgeries on a large class of arborescent knots are
hyperbolic.  Theorem \ref{main} applies to both of these classes of
knots as well as non-alternating, non-arborescent knots and links.

Theorem \ref{main} also gives a nice tool for understanding Dehn
fillings on link complements as well as knot complements.  Classifying
Dehn fillings on links is often a more difficult problem than
classifying fillings on knots, but our arguments apply equally well to
both knots and links.

\subsection{Two proofs, with two notions of length}

By the work of Agol \cite{agol-surg} and Lackenby \cite{lack-surg},
Dehn fillings of a hyperbolic manifold $M$ are hyperbolike whenever
the surgery slopes on $\bdy M$ are ``sufficiently long.'' This term
has two distinct meanings. Agol and Lackenby independently showed that
if the length of each surgery slope on a {\it maximal cusp} of $M$ is
at least $6$, then the surgered manifold is hyperbolike. Lackenby also
showed that the same conclusion holds when the {\it combinatorial
length} of each surgery slope is at least $2\pi$. In this paper, we
use these two points of view to give two different proofs of Theorem
\ref{main}.

Both proofs make use of the same surgery description of the link $K$.
In Section \ref{aug-links}, we show how to start with a prime,
twist-reduced diagram of $K$ and construct a new link, whose
complement in $S^3$ has simple geometric properties. The analysis of
these properties in Section \ref{geom-cusp} leads to estimates of
length on a maximal cusp, yielding our geometric proof of Theorem
\ref{main}. In Section \ref{normal-surf}, we review relevant results
from Lackenby's theory of normal and admissible surfaces in angled
polyhedra, setting up the notion of combinatorial length. We then use
this machinery in Section \ref{our-polyhedra} to give combinatorial
proofs of all three of our main theorems.

We are grateful to Henry Segerman for his helpful suggestions and to
Eric Schoenfeld for providing a template for Figure \ref{fig:rotate}.
Robert Lipshitz has helpfully pointed out that our techniques can be
used to estimate genus.  Above all, we would like to acknowledge the
extended guidance given to us by Steve Kerckhoff. The bulk of these
results were obtained while both authors were his students.

\section{Augmented Links}\label{aug-links}

In this section, we describe how to start with a prime, twist-reduced
projection of a link $K$, construct a flat augmented link $L$, and
subdivide the exterior $S^3 \setminus L$ into two hyperbolic ideal
polyhedra. This construction is originally due to Ian Agol and Dylan
Thurston (see the appendix of \cite{lack-volume}).  We use the
ideal polyhedra to find geometric information about the cusps of the
complement of $L$ in $S^3$.

\subsection{Constructing the augmented link\label{construct-section}}

Let $D(K)$ be a prime, twist-reduced diagram of a link $K \subset
S^3$. As described in the introduction, each twist region in $D(K)$
consists of two strands of $K$ wrapping around each other. For each
twist region $R_i$, add a simple closed curve $C_i$ encircling the
twist region, known as a \emph{crossing circle}.  Let $I$ be the
resulting link.

\begin{define}
  For a link $K \subset S^3$, let the \emph{exterior} $E(K)$ denote
  the complement of an open tubular neighborhood of $K$.
\end{define}

\begin{figure}[b] 
\psfrag{K}{$K$}
\psfrag{I}{$I$}
\psfrag{J}{$J$}
\psfrag{L}{$L$}

\begin{center}
  \includegraphics{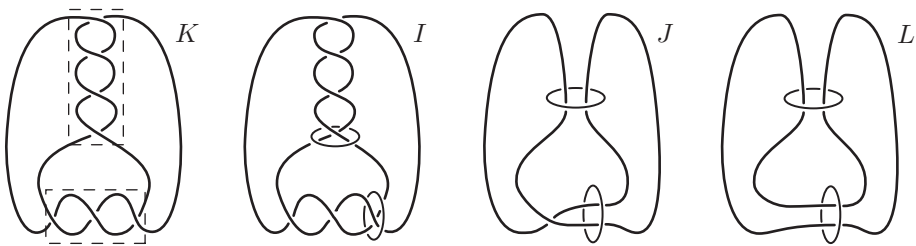}
\end{center}
\caption[Creating the flat augmented link.]{Left to right: The
original knot with twist regions marked; the link $I$ with crossing
circles added; the homeomorphic link $J$; the flat augmented link
$L$.}  \label{fig:aug-construct}
\end{figure}

Note that the manifold $E(I)$ is homeomorphic to the manifold $E(J)$,
where $J$ is the simpler link with all full twists removed at each
twist region of $I$.  We can recover the original link $K$ from $J$ by
performing $1/n_i$ surgery on each $C_i$, $\abs{n_i}$ being the number
of full twists we removed. Furthermore, any Dehn filling of $K$ can be
viewed as a filling of $J$. The results of this paper work by
analyzing the geometry and combinatorics of $S^3 \setminus J$.

In fact, to analyze this geometry, we will make $J$ even simpler by
removing all remaining single crossings from the twist regions. The
resulting link $L$ has two kinds of components: \emph{knot strands}
coming from $K$ that lie flat in the projection plane, and crossing
circles $C_i$ perpendicular to the projection plane. We call $L$ a
\emph{flat augmented link}. If some twist region $R_i$ had an odd
number of crossings, $E(L)$ is no longer homeomorphic to $E(J)$;
indeed, $J$ and $L$ can have a different number of components. We will
address this issue later, in \S \ref{shear-section}. See Figure
\ref{fig:aug-construct} for a visual summary of this construction.

To subdivide $S^3 \setminus L$ into polyhedra, we first slice it along
the projection plane. This divides $S^3$ into two identical $3$-balls.
Since they are identical, we focus our attention on $B_1$, the ball
above the projection plane.  The decomposition of $B_2$ proceeds in
the same way.  Each crossing circle $C_i$ bounds a disk $D_i$, half of
which lies in $B_1$ and borders on three edges in the projection
plane.  We then further slice $B_1$ along each of these half-disks.

\begin{figure}
  \begin{center}
    \includegraphics[width=1.4in]{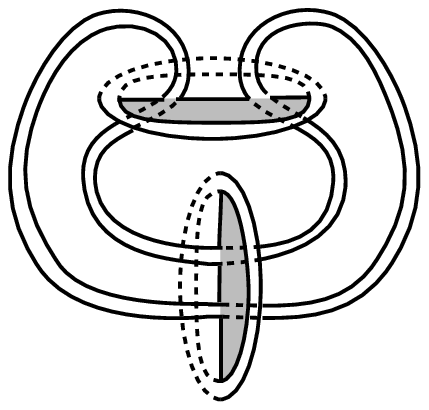}
    \hspace{.1in}
    \includegraphics[width=1.4in]{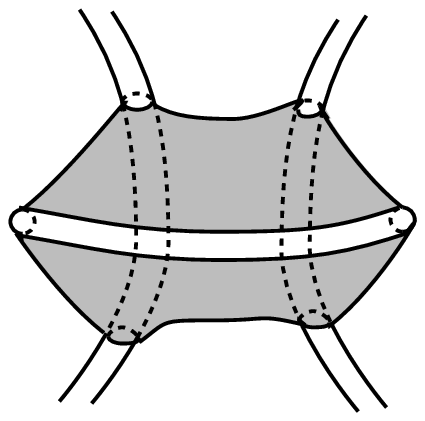}
    \hspace{.1in}
    \includegraphics[width=1.5in]{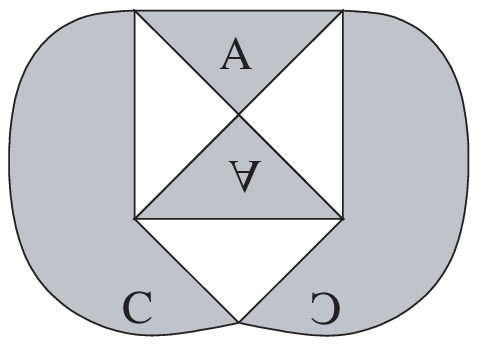}
    \end{center}
    \caption[Decomposing $S^3 \setminus L$ into ideal
    polyhedra.]{Decomposing $S^3 \setminus L$ into ideal polyhedra:
    First slice along the projection plane, then split remaining
    halves of two-punctured disks.  Obtain polygon on right.}
\label{aug-decomp}
\end{figure}

This allows us to pull apart the two sides of each half-disk and
flatten them, creating the planar diagram of a polyhedron. (See Figure
\ref{aug-decomp}.)  This polyhedron will inherit one face from each
region of the projection diagram and one face from each side of disk
$D_i$. We can turn this polyhedron into an \emph{ideal polyhedron}
$P_1$ by collapsing strands of $L$ to ideal vertices. The other ball
$B_2$ becomes an identical ideal polyhedron $P_2$.

$P_1$ and $P_2$ will each have six edges per crossing region, three
from each side of the intersection between $D_i$ and the projection
plane.  At each edge, a face coming from the projection place meets a
face coming from $D_i$. This allows us to two-color the faces in a
convenient fashion: the projection-plane faces will be white and the
crossing-disk faces shaded, as in Figure \ref{aug-decomp}.

To reconstruct $S^3 \setminus L$ from $P_1$ and $P_2$, we first glue
matching shaded faces in each $P_j$, and then glue the two polyhedra
to each other along the white faces. Observe that in this gluing, the
edges become $4$-valent: each borders on two shaded faces (the two
halves of $D_i$) and two white faces in the projection plane. In
$\RR^3$, we can position the crossing disks $D_i$ perpendicular to the
projection plane, creating dihedral angles of $\pi/2$ between adjacent
faces. Conveniently, this feature carries over into hyperbolic
geometry.

\pagebreak

\subsection{The geometry of $E(L)$}

\begin{theorem}[Agol-Thurston \cite{lack-volume}]\label{link-poly-straight}
Let $D(K)$ be a prime, twist-reduced diagram of a link $K$, with at
least two twist regions. Let $L$ be the flat augmented link obtained
from $D(K)$. Then $E(L)$ is hyperbolic. Furthermore, the polyhedra
$P_1$ and $P_2$ decomposing $E(L)$ are convex ideal polyhedra in
$\HH^3$, with totally geodesic faces that meet at right angles.
\end{theorem}

\begin{proof}
Observe that for any flat augmented link $L$, we can always choose a
positive or negative surgery slope $1/n_i$ for each crossing circle
$C_i$ in such a way that filling along these slopes yields an
alternating link $K'$. (This amounts to ensuring that the $n_i$ have
opposite signs for neighboring crossing circles.) Thus every flat
augmented link $L$ is an example of what Adams calls an
\emph{augmented alternating link}. When $L$ has $2$ or more crossing 
circles, and thus $K'$ has $2$ or more twist regions, Menasco's theorem
\cite{menasco-alt} implies $E(K')$ is hyperbolic. Then Adams' result on
augmented alternating links \cite{adams-auglink} implies that every
flat augmented link $L$ is hyperbolic.

Additionally, note that there is an orientation-reversing involution
of $S^3 \setminus L$ preserving $L$ and our ideal polyhedra: namely,
reflection through the projection plane.  Every lift of this
involution to the universal cover $\HH^3$ is a reflection in a totally
geodesic plane. Hence the polyhedra can be made totally geodesic in
$\HH^3$, with the shaded faces meeting the white faces at right
angles.
\end{proof}

\begin{remark}
It is worth noting that the statement and proof of Theorem
\ref{link-poly-straight} do not assume that the original link $K$ is
hyperbolic. When $D(K)$ has at least two twist regions, it follows
from Menasco's theorem \cite{menasco-alt} that the alternating link
$K'$ is hyperbolic; we use this to bootstrap to a hyperbolic structure
on $E(L)$. This will eventually be used to prove that $K$ is
hyperbolic (Theorem \ref{hyp-link}).
\end{remark}

If we intersect $P_1$ and $P_2$ with the compact manifold $E(L)$, each
of their ideal vertices gets truncated into a rectangular
\emph{boundary face} on $\bdy E(L)$. If we keep track of how these
rectangles are glued to one another in the gluing pattern of $P_1$ and
$P_2$, we can construct a picture of the \emph{cusp triangulation} of
each torus of $\bdy E(L)$.

\begin{lemma}\label{flat-cusp-comb}
The cusp tori of $L$ are rectangular. For a crossing circle $C_i$, the
cusp torus is composed of two boundary faces.  For a knot strand $K_j$
lying flat in the projection plane, the cusp torus is a $2 \cross n$
block of boundary faces, where $n$ is the number of twist regions
crossed by $K_j$ (counted with multiplicity).
\end{lemma}

\begin{figure}[b]
\vspace{-10pt}
  \begin{center}
    \includegraphics[width=1.5in]{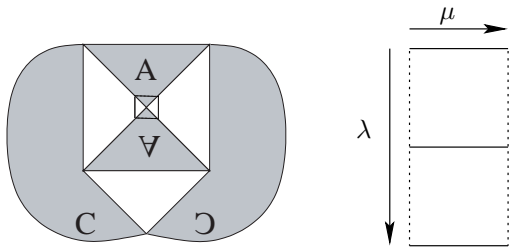}
    \hspace{.25in}
    \input{Fig_rect_small.pstex_t}
  \end{center}

  \caption{Left: Each crossing circle gives rise to one ideal vertex
    in $P_1$.  Right: The cusp diagram corresponding to a crossing
    circle.  Here $\mu$ is a meridian, and $\lambda$ is a longitude.}
    \label{fig:cusp-rect-fig}
\end{figure}

\begin{figure}[ht]
  \begin{center}
    \includegraphics[width=1.5in]{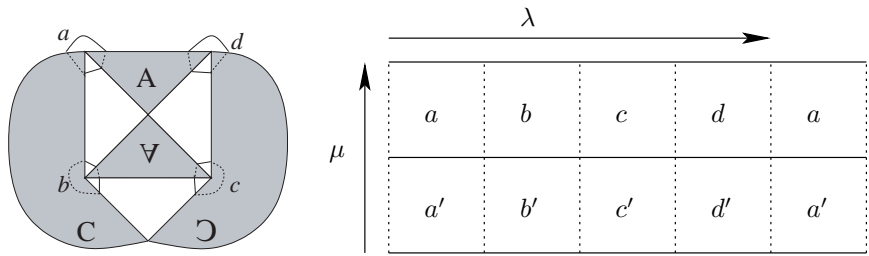}
    \hspace{.1in}
    \input{Fig_rectangle.pstex_t}
  \end{center}
  \caption{The cusp diagram for the knot strand cusp. Solid segments 
correspond to white faces, and dotted segments correspond to shaded faces.}
  \label{fig:knot-strand-rect}
\end{figure}

\begin{proof}
As we saw in the construction of \S \ref{construct-section}, each
crossing circle $C_i \subset L$ becomes an ideal vertex in
$P_1$. Truncate this vertex to get a rectangular boundary face
$F_1$. The shaded faces on opposite sides of this boundary rectangle
are glued to one another in the gluing pattern, since they glue to
give half the disk $D_i$ bounded by $C_i$; thus an arc in $F_1$
connecting the two shaded faces represents a meridian of $C_i$.

The two white faces meeting $F_1$ are glued to corresponding faces of
$P_2$, joining $F_1$ to the boundary rectangle $F_2$. Thus the cusp
torus of $C_i$ is tiled by $F_1$ and $F_2$, with the meridian and
longitude as shown in Figure \ref{fig:cusp-rect-fig}.

For a knot strand $K_j \subset L$, $P_1$ will have one ideal vertex
(hence one boundary rectangle) for each strand of $K_j$ between
adjacent crossing disks $D_i$. (See Figure \ref{aug-decomp}.) These
boundary rectangles are glued end to end along shaded faces coming
from the $D_i$ to complete a longitude of $K_j$. $P_2$ will give rise
to an identical chain of rectangles, glued to the boundary rectangles
of $P_1$ along the white faces of the projection diagram. Thus the
cusp torus of $K_j$ is tiled by a $2 \cross n$ block of rectangles,
where $n$ is the number of intersections between $K_j$ and the
crossing disks $D_i$, hence equal to the number of twist regions that
$K_j$ passes through, counted with multiplicity.  See Figure
\ref{fig:knot-strand-rect}.
\end{proof}

\subsection{Half-twists and surgery slopes}\label{shear-section}
Recall that to construct the flat augmented link $L$ with its nice
polyhedral decomposition, we took three steps, summarized in Figure
\ref{fig:aug-construct}. We added crossing circles to $K$ (obtaining link
$I$); removed a whole number of twists per twist region (obtaining a
homeomorphic link $J$); and then removed any remaining single
crossings. Whereas any Dehn filling of $K$ is a filling of $J$, the
same is no longer true for $L$. Thus to obtain our results, we need to
understand the combinatorics of the link $J$, with the half-twists
re-inserted.

Conveniently, $E(J)$ can still be decomposed into the same polyhedra
$P_1$ and $P_2$, only with a slightly modified gluing pattern. $P_1$
has one shaded face from each side of a crossing disk $D_i$; to
construct $E(L)$, we glued those faces to each other. If instead we
glue those shaded faces of $P_1$ to matching shaded faces of $P_2$, we
effectively insert a half-twist along disk $D_i$ and a single crossing
into the projection diagram of $L$. We can do this wherever $J$ has a
single crossing. In particular, this simple rearrangement means that
we have the following version of Theorem \ref{link-poly-straight}.

\begin{theorem}\label{newlink-poly-straight}
Let $D(K)$ be a prime, twist reduced diagram of a link $K$. Assume
that $D(K)$ has at least two twist regions, with $a_i$ crossings in
twist region $R_i$. Let $J$ be the augmented link constructed in \S
\ref{construct-section}, in which the number of crossings in region
$R_i$ is reduced to $a_i \!\mod 2$. Then
\begin{enumerate}
\item $E(J)$ is hyperbolic, 
\item $E(J)$ is subdivided into convex ideal polyhedra $P_1$ and $P_2$ 
with dihedral angles $\pi/2$,
\item $K$ is the result of Dehn filling each crossing circle $C_i$ of $J$ 
along the surgery slope $1/s_i$, where we removed $2\abs{s_i}$ crossings 
from $R_i$, and
\item every Dehn filling of $K$ is a filling of $J$.
\end{enumerate}
\end{theorem}

\begin{proof}
The first two conclusions follow from Theorem \ref{link-poly-straight}
because $E(J)$ decomposes into the same convex ideal polyhedra as
$E(L)$. The last two conclusions result from the process of
constructing $J$, and are mentioned in \S \ref{construct-section}.
\end{proof}

\begin{figure}[ht]
  \begin{center}
    \includegraphics{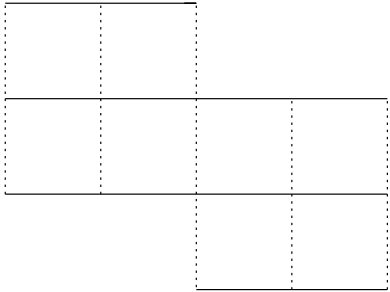}
  \end{center}
  \caption{Cusp view: adding a half twist to a flat augmented link shifts the
  gluing along the shaded faces.}
  \label{fig:cusp_twist}
\end{figure}

In the cusp diagrams of $\bdy E(J)$, each half-twist in the transition
from $L$ to $J$ will shift the gluing by one step along the shaded
faces coming from $D_i$, as illustrated in Figure
\ref{fig:cusp_twist}.  Thus the neat rectangular pattern of Lemma
\ref{flat-cusp-comb} no longer holds. However, we can still make
convenient statements about the universal cover $\tild{T}$ of
each cusp torus of $S^3 \setminus J$.

\begin{define}\label{lattice-basis}
Let $T$ be a cusp torus of $\bdy E(J)$, with universal cover $\tild{T}
= \RR^2$. Then $\tild{T}$ contains a rectangular lattice coming from
white and shaded faces of $P_1$ and $P_2$. We construct a basis
$\langle \s, \w \rangle$ of this $\ZZ^2$ lattice by letting $\s$ be a
step parallel to a shaded face and $\w$ be a step parallel to a white
face.
\end{define}

\begin{lemma}\label{new-cusp-comb}
Let $T$ be a cusp torus of $\bdy E(J)$ and let $\langle \s, \w
\rangle$ be the basis for the lattice on $\tild{T}$. In this basis,
the fundamental domain of $T$ appears as follows:

\begin{enumerate}
\item If $T$ comes from a crossing circle without a half-twist, it has 
meridian $\w$ and longitude $2\s$.
\item If $T$ comes from a crossing circle with a half-twist, it has meridian 
$\w \pm \s$ (depending on the direction of the twist) and longitude $2\s$.
\item If $T$ comes from a component $K_j$ of the original link $K$,
it has meridian $2\s$ and longitude $n\w+k\s$, where $K_j$ runs
through $n$ twist regions with multiplicity and $k$ is an undetermined
integer.
\end{enumerate}
\end{lemma}

\begin{proof}
If $J$ does not contain any half-twists, this is a restatement of
Lemma \ref{flat-cusp-comb}. (See Figures \ref{fig:cusp-rect-fig} and
\ref{fig:knot-strand-rect}.)  Each half-twist along the crossing circle
$C_i$ shears the meridian of $C_i$ by $\s$, as described above.  It
also shears by $\s$ the cusp of every component of the original link
$K$ passing through the crossing disk $D_i$.

Thus if $K_j$ passes through a half-twist $m$ times, the projection of
the curve $n\w+m\s$ to $T$ will be \emph{some} longitude of $K_i$, in
the sense of completing a basis of $\pi_1(T)$ along with $\mu=2\s$.
The \emph{true} longitude, in the sense of having linking number $0$
with $K_j$, is then some curve of the form $n\w+k\s$ for some integer $k$.
\end{proof}

The basis $\langle \s, \w \rangle$ also allows us to make precise
statements about the surgery curves on $\bdy E(J)$ that correspond to
non-trivial surgeries on $K$.

\begin{theorem}\label{surg-curve-in-basis}
Let $K = \cup_{j=1}^m K_j$ be a link in $S^3$ with a prime, twist
reduced diagram $D(K)$. Suppose that $D(K)$ contains twist regions
$R_1, \ldots, R_n$ ($n \geq 2$) and that twist region $R_i$ contains
$a_i$ crossings. For each component $K_j$, let $n_j$ be the number of
twist regions crossed by $K_j$, counted with multiplicity; and let
$s_j$ be a non-trivial surgery slope on $K_j$.

With this notation, the surgery on $S^3 \setminus K$ along slopes
$s_1, \ldots, s_m$ can be represented as a surgery on $J$ as follows:

\begin{enumerate}
\item On the (mostly) planar component of $J$ corresponding to $K_j$, 
the surgery curve is $p_j n_j \w + q_j \s$, for some integers $p_j \neq
0$ and $q_j$.
\item On the crossing circle $C_i$, the surgery curve is $\w \pm a_i \s$.
\end{enumerate}
\end{theorem}

\begin{proof}
By Lemma \ref{new-cusp-comb}, $K_j$ has meridian $2\s$ and a longitude
of the form $n_j \w + k_j \s$. Since $s_j$ is a non-trivial surgery
slope, it must cover at least one longitude. In particular, the number
of steps that a curve representing $s_j$ takes along the white faces
is a nonzero multiple of $n_j$.

To prove conclusion $(2)$, suppose first that $a_i$ is even, so $J$
has no half-twist at $C_i$. By Lemma \ref{new-cusp-comb}, $C_i$ has
longitude $2\s$ and meridian $\w$. By Theorem
\ref{newlink-poly-straight}, the surgery curve on $C_i$ traverses
$a_i/2$ longitudes and one meridian, proving the result.

Now, suppose that $a_i$ is odd.  Then in the construction of $J$, we
have removed $2b_i = a_i-1$ crossings; the remaining half-twist of $J$
at $C_i$ goes in the same direction as the twists of $K$.  By Lemma
\ref{new-cusp-comb}, $C_i$ has longitude $2\s$ and meridian 
$\w + \sigma_i \s$, for some $\sigma_i = \pm 1$.  By Theorem
\ref{newlink-poly-straight}, the surgery curve traverses $\sigma_i
b_i$ longitudes (with the same $\sigma_i$) and one meridian.  Thus, in
the basis of $\langle \s, \w \rangle$, the surgery curve is
{\setlength\arraycolsep{3pt}
\setlength{\belowdisplayskip}{-1ex}
\begin{eqnarray*}
\mu + \sigma_i b_i \lambda 
&=& (\w + \sigma_i \s) + \sigma_i b_i (2\s) \\
&=& \w + \sigma_i (1 + 2b_i) \s \\
&=& \w + \sigma_i a_i \s \, .
\end{eqnarray*}}
\end{proof}

\section{Geometric Cusp Estimates}\label{geom-cusp}
In Section \ref{aug-links}, we showed that each cusp of $E(J)$
contains a rectangular lattice generated by $\s$ and $\w$, and in
Theorem \ref{surg-curve-in-basis} we expressed the surgery curves on
$\bdy E(J)$ explicitly in terms of $\s$ and $\w$.  In this section, we
will use the geometry of the polyhedra $P_1$ and $P_2$ to come up with
lower bounds for the lengths of $\s$ and $\w$ on a maximal cusp. This
will allows us to estimate the lengths of surgery curves. By combining
these estimates with Agol and Lackenby's $6$-Theorem (Theorem
\ref{6-theorem}), we will obtain a geometric proof of Theorem \ref{main}.

\subsection{Length on a maximal cusp}\label{max-section}
In this paper, we measure the length of curves on a cusp of $E(J)$ in
two distinct ways: geometric and combinatorial.  The geometric
measurements of this section come from the hyperbolic metric.  A
closed curve isotopic to a cusp torus has many representatives in
$E(J)$, whose lengths shrink to $0$ as the representative curves
approach the cusp.  To obtain a meaningful definition of length, we
will consider curves on a horospherical torus bounding a maximal
neighborhood of a cusp.

For a manifold with just one cusp, such as a knot complement, we
obtain the maximal horoball neighborhood by expanding a horoball about
the single cusp until it becomes tangent to itself.  In a manifold
with multiple cusps, such as $E(J)$, the size of the maximal cusp
depends on the order in which we expand horoballs about the cusps, for
a horoball might become tangent to an expanded horoball about a
different cusp before it becomes tangent to itself.

\begin{define}\label{geom-length-def}
Let $M$ be a $3$-manifold with boundary consisting of tori, whose
interior has a complete hyperbolic structure. Fix a {\it cusp
neighborhood} $U$, consisting of disjoint horoball neighborhoods of
the cusps of $M$. Then any closed curve $\gamma \subset \bdy M$ can be
assigned a unique {\it geometric length} $\lgeom(\gamma)$, defined as
the shortest length of a curve on $\bdy U$ isotopic to $\gamma$.

The subscript in $\lgeom(\gamma)$ serves to distinguish geometric
length from the combinatorial length $\lcomb(\gamma)$ used in Sections
\ref{normal-surf} and \ref{our-polyhedra}. When the meaning is clear,
we will simply use $\ell(\gamma)$.
\end{define}

To rule out exceptional surgeries on $M$, it helps to choose the cusp
neighborhood $U$ to be maximal, to make the surgery curves as long as
possible.  Agol \cite{agol-surg} and Lackenby \cite{lack-surg} have
independently proved the following surgery theorem:

\begin{theorem}[\textbf{6-Theorem}]\label{6-theorem}
Let $M$ be a hyperbolic $3$-manifold with boundary consisting of
tori. Let $s_1, \ldots, s_n$ be surgery slopes on $\bdy M$, with one
$s_i$ on each torus. Suppose that there are disjoint horoball
neighborhoods of the cusps of $M$, such that $\lgeom(s_i) > 6$ for all
$i$. Then the manifold obtained by Dehn filling $M$ along the slopes
$s_1, \ldots, s_n$ is hyperbolike.
\end{theorem}

In \S \ref{horosphere-section}, we will give explicit instructions for
expanding the horoball neighborhoods about the cusps of $E(J)$ that
produce favorable estimates for the length of surgery curves.  These
estimates will rely on lower bounds for $\ell(\s)$ and $\ell(\w)$.

It should be noted that the lengths of the meridians of augmented
links were found independently by Eric Schoenfeld in his undergraduate
thesis \cite{schoenfeld}.

\subsection{Horosphere packing in $\HH^3$}\label{horosphere-section}

Recall from Section \ref{aug-links} that each ideal vertex of $P_1$
and $P_2$ gives rise to a {\it boundary rectangle} on a cusp of
$E(J)$. In the geometry of the universal cover $\HH^3$, the boundary
rectangle can be seen as the intersection of $P_i$ with a horosphere.
A side parallel to $\s$ is the intersection of a horosphere with a
shaded face, and a side parallel to $\w$ is the intersection with a
white face. It turns out that boundary rectangles and horospheres are
easiest to visualize in the upper half-space model of $\HH^3$.

\begin{notation}
We will parameterize the upper half-space model of $\HH^3$ with coordinates
$(z,h)$, where $z \in \CC$ and $h \in \RR^+$. In this model, the sphere at
infinity $S^2_\infty$ can be identified with the Riemann sphere $\CC \cup 
\infty$.
\end{notation}

We can apply an isometry of $\HH^3$ so that the point at infinity of
$\HH^3$ projects to the cusp under consideration in $E(J)$.  Then a
horosphere about that cusp lifts to a horizontal plane 
at height $h$. In the metric on the upper half-space model, hyperbolic
length corresponds to $\frac{1}{h}$ times the Euclidean
length. Thus we can bound the hyperbolic lengths of $\s$ and $\w$
using Euclidean measurements in upper half-space.

Every procedure for expanding the cusps will lower the horizontal
horosphere until it becomes tangent to another horosphere. This
abutting horosphere will look like a Euclidean sphere tangent to
$S^2_\infty$ at some point of $\CC$.

\begin{define}\label{def:horosphere-center}
Let $H$ be a horosphere in the upper half-space model of $\HH^3$. If
$H$ is a Euclidean sphere, call the point of tangency in $\CC$ the
{\it center} of $H$.  If $H$ is a horizontal Euclidean plane, we say
that $H$ is {\it centered at $\infty$}.
\end{define}

We will normalize our horoball packing in upper half-space by placing 
ideal triangles into {\it standard position}.

\begin{define}
Let $T \subset \HH^3$ be an ideal triangle. We will say that $T$ is in
{\it standard position} in the upper half-space model if its ideal
vertices lie at $0$, $1$, and $\infty$. Note that any ideal triangle
can be placed into standard position by an isometry.
\end{define} 

Recall that each shaded face of the polyhedra $P_1$ and $P_2$ is an
ideal triangle. Thus we can place $P_i$ in the upper half-space model
of $\HH^3$ so that a shaded face is in standard position.  When we do
so, the polyhedron $P_i$ is lifted to lie entirely over a boundary
rectangle with corners at $0$, $1$, $0+ir$ and $1+ir$ for some real
number $r$.

\begin{lemma}\label{lemma:off-center-horo}
Arrange $P_i$ in $\HH^3$ so that some shaded side has vertices at $0$,
$1$, and $\infty$, so $P_i$ lies over a boundary rectangle with
corners at $0$, $1$, $0+ir$ and $1+ir$.  Follow any procedure for
expanding the cusps of $E(J)$ to horoballs with disjoint interiors,
and let $H$ be a horosphere centered at a point of the rectangle. If
the center of $H$ does not lie on a white side of the rectangle, the
diameter of $H$ is at most $1$.
\end{lemma}

\begin{proof}
Since there is a shaded face of $P_i$ with vertices $0$ and $1$, there
must be a white face of $P_i$ containing vertices $0$, $ir$ and
$\infty$, and another containing $1$, $1+ir$, and $\infty$.  Recall
that reflection through the white faces of $E(J)$ is an involution of
the manifold, corresponding to a reflection in the projection plane of
Figure \ref{aug-decomp}.  This involution permutes the horospheres
covering cusps of $E(J)$, and thus takes $H$ either to itself or to a
disjoint horosphere $H'$.  If the center of $H$ lies on a white side
of the rectangle, reflection in the plane above that side will fix
$H$.  Otherwise, $H$ must be disjoint from its reflection.  Then,
since the boundary rectangle has length $1$, the diameter of $H$ can
be at most $1$.
\end{proof}

For horospheres centered on white sides of boundary rectangles, we
will also prove that the diameter is at most $1$. In order to do so,
we will need to give specific instructions for expanding the cusps.

We would like to expand the cusps of $E(J)$ to a halfway point along
each edge of the polyhedra. It turns out that even though
an edge is infinitely long, there is a natural way to define
its midpoint.

\begin{define}\label{def:midpoint}
Let $T \subset \HH^3$ be an ideal triangle. For each edge $e$ of $T$,
define the {\it midpoint} to be the point $m \in e$ such that the
geodesic from $m$ to the opposite vertex is perpendicular to $e$.
(This point is unique, for otherwise we would have a triangle with one
ideal vertex and two right angles.)
\end{define}

Now, each edge $e$ in the polyhedral decomposition of $E(J)$ borders
on two shaded faces, $S_1$ and $S_2$, with each $S_i \subset P_i$. (See
Figure \ref{aug-decomp}.) It is easy to check that the two definitions
of the midpoint of $e$, coming from $S_1$ and $S_2$, coincide. This is
because $P_1$ and $P_2$ are symmetric by a reflection in the white
faces of the projection diagram. The reflection preserves angles, so
it also preserves the midpoint of $e$. Thus we have a well-defined
midpoint of each edge of $E(J)$.

\begin{figure}[ht]
  \begin{center} 
  \input{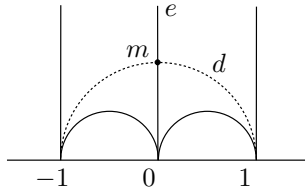} 
  \end{center}
  \caption{When a triangle is in standard position, the midpoint of a
  vertical edge lies at height $1$.}  
  \label{fig:reflect}
\end{figure}

When two ideal triangles $T_1$ and $T_2$ are symmetric across an edge
$e$, we have an alternate way of seeing the midpoint $m$ in $e$.
Namely, the two ideal triangles glue up to form an ideal
quadrilateral.  One diagonal is $e$, and the other diagonal $d$
intersects $e$ at $m$. (See Figure \ref{fig:reflect}.)  By this
approach, we see that when $T_1$ is in standard position and $e$ is
vertical, $m$ is at Euclidean height $1$.

We are now ready to expand all the cusps of $E(J)$.  Pick an order for
the cusps, $K_1$, $K_2$, $\ldots$, $K_r$, and expand one $K_i$ at a
time, starting with $K_1$.  Continue expanding the horoball
neighborhood of $K_i$ until it either meets another horoball, or meets
the midpoint of some edge into $K_i$.

\begin{lemma}\label{lemma:white-side-horo}
Arrange $P_1$ or $P_2$ in $\HH^3$ so that some shaded face is in
standard position, and let $H$ be a horoball centered on a white side
of a boundary rectangle of $P_i$.  If the interior of $H$ does not
contain the midpoint of any edge of $E(J)$, the diameter of $H$ can
be at most $1$. 
\end{lemma}

\begin{proof}

The center of $H$ on $S^2_\infty$ is an ideal vertex of $P_i$, and so
at this ideal vertex, two white faces and two shaded faces meet and
intersect $H$ in a rectangle.  Recall that these faces are all totally
geodesic, by Theorem \ref{newlink-poly-straight}.  Since $H$ is
centered on a white side of a boundary rectangle, we know one of these
white faces is actually vertical.  That is, one white face $V$ meeting
the center of $H$ lies in a vertical plane in $\HH^3$, bordered by a
line $\ell \subset \CC$.  The other white face, $W$, must also lie in
a geodesic plane $P$ in $\HH^3$.  $P$ cannot also be vertical (since
$H$ is not centered at $\infty$), so it looks like a Euclidean
half-sphere, tangent to the vertical plane containing $V$.  Thus the
boundary of $P$ at infinity is a circle $C$.  (See Figure
\ref{fig:rotate}.)

Consider the white face $W$.  Since white faces of $P_i$ only meet
white faces at ideal vertices, other white faces meeting $W$ will lie
on geodesic planes in $\HH^3$ tangent to $P$. Thus they extend to give
circles tangent to the circle $C$.  In particular, these circles have
disjoint interiors, and so the interior of $C$ must be disjoint from
the two vertical planes containing the white sides of the boundary
rectangle about the vertex of $P_i$ at $\infty$.  But because we put a
shaded side in standard position, these vertical planes are of
Euclidean distance $1$ apart.  Thus the diameter of $C$ can be no more
than $1$.

\begin{figure}
  \begin{center}
\psfrag{H}{$H$}
\psfrag{C}{$C$}
\psfrag{l}{$\ell$}
\psfrag{g}{$\gamma$}
\psfrag{e}{$e$}
\psfrag{z}{$z$}
\psfrag{m}{$m$}
\includegraphics{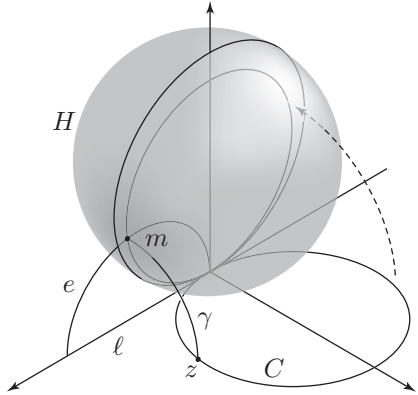}
  \end{center}
  \caption{If the diameter of $H$ is greater than $1$, a $90^\circ$ 
  rotation about the side of the boundary rectangle shows that the 
  midpoint of $e$ lies inside $H$.}
  \label{fig:rotate}
\end{figure}

Now consider the two shaded faces meeting at the center of the
horosphere $H$.  At most one of these shaded faces can lie in a
vertical plane.  Let $S$ be one of these shaded faces that doesn't lie
in a vertical plane.  Let $e$ be the edge between $S$ and the vertical
white face $V$, and let $z$ be the vertex of $S$ opposite $e$.  Note
that $z$ is an endpoint of the edge given by the intersection of $S$
with the white face $W$. (The other endpoint is the center of $H$.)
Thus $z$ lies on the circle $C$.

Consider the midpoint $m$ of $e$. By Definition \ref{def:midpoint}, the 
geodesic $\gamma$ from $z$ to $m$ meets $e$ at a right angle.
But $\gamma$ must lie in $S$, which by Theorem \ref{newlink-poly-straight}
meets the face $V$ at a right angle, and hence $\gamma$ is normal to
the entire vertical plane containing $V$.  Thus $\gamma$ is a
Euclidean quarter circle, centered on the line $\ell$ in $\CC$ that
lies under the vertical white face $V$.  Then note that a $90^\circ$
rotation about $\ell$ takes the point $z$ to $m$.  But this same
$90^\circ$ rotation about $\ell$ will take the circle $C$ to the
vertical plane over $\ell$.  If the diameter of $H$ is greater than
$1$, the circle $C$, of diameter at most $1$, will be contained inside
$H$.  Since $z$ is on $C$, in this case the rotated point $m$ will be
contained inside $H$.  (See Figure \ref{fig:rotate}.)  Hence if we do
not allow $H$ to contain $m$, the diameter of $H$ can be at most $1$.
\end{proof}

\begin{theorem}\label{thm:tangent-at-midpoint}
Expand all the cusps of $E(J)$ as above.  Then the midpoint of every
edge of $E(J)$ will lie at the point of tangency of two horospheres.
\end{theorem}

\begin{proof}
First, we would like to show that an expanding horoball about a given
cusp of $E(J)$ will simultaneously meet the midpoints of all the edges
into that cusp.  To that end, consider the horoball $H_\infty$ about
$\infty$, normalized so that a shaded face $S$ is in standard
position.  $S$ forms a side of a boundary rectangle of Euclidean width
$1$.  The opposite side of this rectangle must be another shaded side
of Euclidean width $1$.  Continuing in the $\w$ direction, we see
there is an infinite strip consisting of boundary rectangles lined end
to end, and each shaded side in this strip has Euclidean width $1$.

We can now reflect this infinite strip by an involution in the white
faces, obtaining an infinite strip whose sides are separated by $2\s$.
By Lemma \ref{new-cusp-comb}, translation by $2\s$ is an element of
the covering transformation group for any cusp of $E(J)$.  Thus every
shaded face intersecting $H_\infty$ has Euclidean width $1$.
Consequently, the midpoint of every vertical edge is at height $h=1$.
(See Figure \ref{fig:reflect}.)

As we expand the horoball $H_\infty$ about $\infty$, the expansion
stops before height $1$ only if $H_\infty$ becomes tangent to another
horoball $H$ of diameter greater than $1$.  But by Lemmas
\ref{lemma:off-center-horo} and \ref{lemma:white-side-horo}, all
horoballs obtained with our expansion instructions have diameter at
most $1$.  Thus $H_\infty$ can be expanded until it reaches height
$h=1$ and meets the midpoint of every vertical edge.  By symmetry, the
horosphere centered at the other endpoint of a vertical edge also
meets the midpoint of that edge.
\end{proof}

\begin{corollary}\label{cor:s=1}
Expand the cusps of $E(J)$ by the above procedure. Then in any boundary
rectangle, $\ell(\s)=1$ and $\ell(\w) \geq 1$.
\end{corollary}

\begin{proof}
Consider a boundary rectangle on a maximal cusp of $E(J)$, positioned
so that a shaded face $S$ adjacent to this rectangle is in standard
position.  By Theorem \ref{thm:tangent-at-midpoint}, the horosphere
$H_\infty$ about $\infty$ meets the the midpoints of the vertical
edges of $S$.  In standard position, these midpoints lie at $h=1$,
where Euclidean lengths correspond to hyperbolic lengths.  Thus
$\ell(\s)$ is the Euclidean width of $S$, namely $1$.

Theorem \ref{thm:tangent-at-midpoint} also implies that above every
corner of a boundary rectangle, horosphere $H_\infty$ is tangent to a
horosphere of diameter $1$. Since these horospheres are disjoint, we
can conclude that $\ell(\w) \geq 1$.
\end{proof}

\subsection{Surgery consequences}

By combining Corollary \ref{cor:s=1} with Theorem
\ref{surg-curve-in-basis}, we can compute explicit lower bounds for
the lengths of surgery curves.

\begin{theorem}\label{thm:curve-lengths}
For each cusp $K_i$ of $E(J)$, pick a surgery slope $s_i$ that
represents a non-trivial filling of the original link $K$. If $K_i$ is
a knot strand cusp, let $n_i$ be the number of twist regions visited by
the corresponding strand of $K$; if $K_i$ is a crossing circle cusp,
let $n_i$ be the number of crossings in the corresponding twist region.
Then
\begin{enumerate}
\item on a knot strand cusp, $\lgeom(s_i) \geq n_i$, and
\item on a crossing circle cusp, $\lgeom(s_i) \geq \sqrt{n_i^2 + 1}$.
\end{enumerate}
\end{theorem}

\begin{proof}
For part (1), Theorem \ref{surg-curve-in-basis} implies that any
non-trivial surgery curve on $K_i$ is of the form $p_i n_i\w+q_i \s$ for
integers $p_i \neq 0$ and $q_i$.  Then since $\ell(\s)=1$ and
$\ell(\w) \geq 1$, and $\s$ and $\w$ are perpendicular, any such curve
will have length at least $\sqrt{p_i^2 n_i^2+q_i^2}$.  This is minimal
when $q_i=0$ and $p_i = \pm 1$.  In this case, $\ell(s_i) \geq n_i$.

For part (2), Theorem \ref{surg-curve-in-basis} implies the surgery
curve is $\w \pm n_i\s$, which has length at least $\sqrt{n_i^2 +1}$.
\end{proof}

We are now ready to give our geometric proof of Theorem \ref{main},
which we restate.

\medskip

\noindent\textbf{Theorem \ref{main}.} \emph{
Let $K$ be a link in $S^3$ with a prime, twist-reduced diagram
$D(K)$. Suppose that every twist region of $D(K)$ contains at least
$6$ crossings and each component of $K$ passes through at least $7$
twist regions (counted with multiplicity). Then every non-trivial Dehn
filling of all the components of $K$ is hyperbolike.
}

\begin{proof}
Since each knot strand cusp crosses at least $7$ twist regions, Theorem
\ref{thm:curve-lengths} says that the surgery curve on that cusp has 
length at least $7$. Since each twist region contains at least $6$
crossings, the surgery curve on the corresponding crossing circle has
length at least $\sqrt{6^2+1} > 6$. Thus the surgery curve on every
component of $\bdy E(J)$ has length greater than $6$. Therefore, by
the 6-Theorem, the surgered manifold is hyperbolike.
\end{proof}

\section{Angled Polyhedra and Normal Surfaces}\label{normal-surf}

The next two sections give a combinatorial proof of Theorems
\ref{hyp-link} and \ref{main}, by using a combinatorial notion of
length.  In this proof, we make use of a number of results from the
theory of normal and admissible surfaces in angled polyhedra, much of
it developed by Marc Lackenby \cite[Section 4]{lack-surg}.  Lackenby
worked with the dual structure of angled spines.  We find it more
convenient to work with polyhedra, so in this section we will
translate his definitions and theorems into polyhedral language.

\subsection{Normal and admissible surfaces}

\begin{define}
For the purposes of this paper, a {\it polyhedron} is a $3$-ball $P$
with a specified graph $\Gamma$ embedded in $\bdy P$, such that
\begin{enumerate}
\item each vertex of $\Gamma$ has valence at least $3$,
\item each edge of $\Gamma$ has ends on distinct vertices, and
\item each region of $\bdy P \setminus \Gamma$ is bounded by at least $3$ 
edges.
\end{enumerate}
$P$ inherits vertices and edges from $\Gamma$, and the {\it faces} of
$P$ are regions of $\bdy P \setminus \Gamma$. An {\it ideal
polyhedron} is a polyhedron with the vertices removed.
\end{define}

\begin{remark}
This definition of an ideal polyhedron is actually slightly stronger
than Lackenby's dual definition of a thickened spine, in that
condition ($3$) prohibits our polyhedra from having bigon faces. This
stronger definition is sufficient for our purposes (certainly, the
ideal polyhedra $P_1$ and $P_2$ constructed in Section \ref{aug-links}
have no bigon faces), and allows for stronger statements of some
results.
\end{remark}

Let $M$ be a manifold subdivided into ideal polyhedra.  To see how the
ideal vertices fit together to tile $\bdy M$, we truncate all the
ideal vertices.  This gives new polyhedra, with two kinds of faces:
{\it interior faces} that are truncated copies of the original faces,
and {\it boundary faces} that come from the truncated vertices.  We
also obtain two kinds of edges: {\it interior edges} that come from
the original truncated edges, and {\it boundary edges} along the
boundary faces.

In order to define a combinatorial length for a curve on $\bdy M$, we
actually need to consider surfaces inside the manifold with that curve
as boundary.  Thus we review some results from the theory of normal
and admissible surfaces.

\begin{figure}
\begin{center}
\psfrag{(a)}{(a)}
\psfrag{(b)}{(b)}
\includegraphics{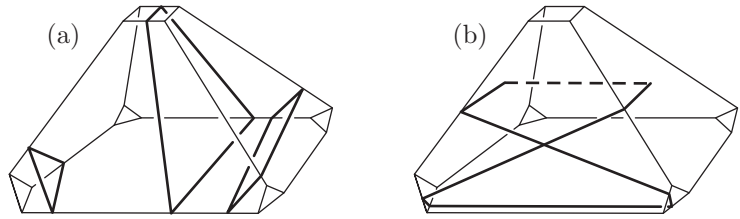}
\caption{(a) Normal disks in a truncated polyhedron. (b) An admissible disk.}
\label{normal-disk}
\end{center}
\end{figure}

Let $(F, \bdy F) \subset (M, \bdy M)$ be an embedded essential surface
(a sphere not bounding a ball, or an incompressible,
boundary-incompressible surface).  The theory of normal surfaces,
originally developed by Haken \cite{haken-normal} and generalized and
expanded in many directions, says that $F$ can be isotoped until its
intersections with the polyhedra (or handles) have a particularly nice
form.  Specifically, we can get $F$ to intersect each polyhedron in a
collection of disjoint, embedded disks, with each disk positioned so
that its boundary curve $\gamma$ has several nice properties:

\begin{define}\label{normal-def}
Let $P$ be a truncated ideal polyhedron.  A simple closed curve
$\gamma \subset \bdy P$ is called {\it normal} if
\begin{enumerate}
\item $\gamma$ is transverse to the edges of $P$,
\item no arc of $\gamma$ in a face of $P$ has endpoints on the same edge, 
or on an interior edge and an adjacent boundary edge,
\item $\gamma$ doesn't lie entirely in a face of $P$,
\item $\gamma$ intersects each edge at most once, and
\item $\gamma$ intersects each boundary face at most once.
\end{enumerate}

The disk in $P$ bounded by a normal curve $\gamma$ is called a {\it
normal disk}. See Figure \ref{normal-disk}(a) for several examples.
\end{define}

\begin{notation}
To avoid confusion with longer arcs on $\bdy M$, we will refer to the
arcs of intersection between a normal curve $\gamma$ and the faces of
$P$ as {\it segments}. (Thus an arc can consist of many segments.) 
Segments of $\gamma$ lying in interior faces of $P$ will be called {\it
interior segments}, and the segments lying in boundary faces will be
called {\it boundary segments}.
\end{notation}

In order to prove word-hyperbolicity, we actually need to work with a
more general class of surfaces that cannot be normalized. These
surfaces may not be embedded, and may even have boundary components in
the interior of $M$.

\begin{define}\label{adm-def}
Let $P$ be a truncated ideal polyhedron. An immersed disk $D \subset P$ 
is called {\it admissible} if
\begin{enumerate}
\item $\bdy D \setminus \bdy P$ is a (possibly empty) collection of 
embedded arcs with endpoints inside interior faces of $P$,
\item $\bdy D \cap \bdy P$ is an immersed closed curve or an immersed 
collection of arcs,
\item if $\bdy D \cap \bdy P$ is a closed curve, it satisfies
conditions $(1) - (3)$ of Definition \ref{normal-def} of a normal curve,
\item each arc component of $\bdy D \cap \bdy P$ satisfies conditions 
$(1) - (2)$ of Definition \ref{normal-def}, and
\item each segment of $\bdy D$ in a face of $P$ is embedded.
\end{enumerate}
\noindent 
An example is shown in Figure \ref{normal-disk}(b). 
We call an immersed surface $F \subset M$ an {\it admissible surface}
if it intersects each polyhedron in a collection of admissible disks.
\end{define}

\subsection{Angle structures and combinatorial area}

The theory of normal surfaces becomes much more powerful if one has
information about the dihedral angles of ideal polyhedra.

\begin{define}\label{angled-def}
Let $M$ be a manifold with boundary. An {\it angled polyhedral
decomposition} of $M$ is a subdivision of $M \setminus \bdy M$ into
ideal polyhedra, glued along their interior faces. Each interior edge
of each {\it angled polyhedron} is assigned an internal angle
$\alpha_i \in (0, \pi)$ and an external angle $\epsilon_i = \pi -
\alpha_i$, such that
\begin{enumerate}
\item around each edge of $M$, $\sum \alpha_i = 2 \pi$, and
\item in each polyhedron, for a normal curve $\gamma$ that intersects only 
interior edges, $\sum_{\gamma} \epsilon_i \geq 2\pi$, with equality if and 
only if $\gamma$ encircles a vertex.
\end{enumerate}
\end{define}

Angle structures on a polyhedral decomposition of $M$ allow us to
define the combinatorial area of a surface.

\begin{define}\label{area-def}
Let $D \subset P$ be an admissible disk in an angled polyhedron, with
the boundary faces of $P$ lying on $\bdy M$. Let $E_1, \ldots, E_n$ be
the interior edges crossed by $\bdy D$ (counted with multiplicity),
and let $\epsilon_1, \ldots, \epsilon_n$ be the corresponding external
angles. Then define the combinatorial area of $D$ to be
$$a(D) = \sum_{i=1}^n \epsilon_i + \pi \abs{\bdy D \cap \bdy M} 
- 2 \pi + 3 \pi \abs{\bdy D \setminus \bdy P} \, .$$

For an admissible surface $F \subset M$, $a(F)$ is defined by summing the 
areas of its admissible disks.
\end{define}

For disks with $\bdy D \subset \bdy P$, this definition matches the
formula for hyperbolic area. For a polygon $T \subset \HH^2$ with
external angles $\epsilon_i$, $a(T) = \sum \epsilon_i - 2\pi$. (See,
for example, Corollary 2.4.15 of \cite{thur-book}.)  Ideal vertices
have internal angle $0$ and thus add $\pi$ to the area, just as each
component of $\bdy D \cap \bdy M$ adds $\pi$ to combinatorial area. As
for the coefficient $3\pi$ per component of $\bdy D \setminus \bdy P$,
it was chosen by Lackenby to make the combinatorial area of $D$
automatically positive whenever $\bdy D$ passes through the interior
of $P$.

\begin{figure}
\begin{center}
\psfrag{(a)}{(a)}
\psfrag{(b)}{(b)}
\includegraphics{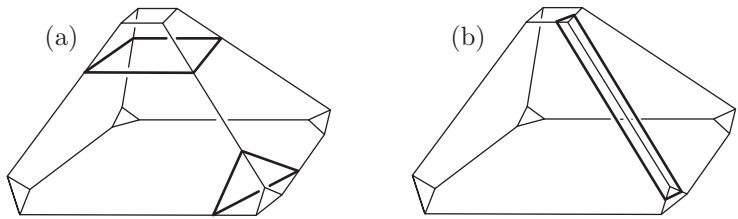}
\caption{(a) Vertex links. (b) A boundary bigon.}
\label{zero-area-fig}
\end{center}
\end{figure}

In fact, there are only two types of admissible disks whose area is
$0$; both of them happen to be normal. The first is a {\it vertex
link} cutting off a boundary face; its area is $0$ by Definition
\ref{angled-def}. The second is a {\it boundary bigon} cutting off an
interior edge; it has area $0$ because its boundary curve only picks
up area from two boundary faces. They are shown in Figure
\ref{zero-area-fig}.

\begin{lemma}\label{pos-area}
Let $D \subset P$ be an admissible disk in an angled polyhedron. If
$D$ is not a vertex link or a boundary bigon, $a(D) > 0$.
\end{lemma}

\begin{proof}
If $D$ is normal, Lackenby \cite[Lemma 4]{lack-bundle} proves that
$a(D) > 0$ unless $D$ is a vertex link or a bigon. If $D$ is not
normal, Lackenby \cite[Lemma 4.2]{lack-surg} proves that $a(D) >
0$. In both cases, the proofs rely on condition $(2)$ of Definition
\ref{angled-def} of an angled polyhedron and the observation that if
$\bdy D$ self-intersects or crosses an edge multiple times, the area
can actually be reduced by surgering the disk.
\end{proof}

\begin{remark}
Lemma \ref{pos-area} is one place where our definition of a
polyhedron, which is stronger than Lackenby's definition of an angled
spine because it rules out bigon faces, becomes convenient.  Bigon
faces of the polyhedra allow normal disks other than boundary bigons
or vertex links to have zero area \cite[Lemma 4.2]{lack-surg}; in our
scenario, every other admissible disk has strictly positive area.
\end{remark}

The analogy between hyperbolic area and combinatorial area extends to
the following combinatorial version of the Gauss-Bonnet theorem
\cite[Proposition 4.3]{lack-surg}.

\begin{prop}[\textbf{Gauss-Bonnet Theorem}]\label{gauss-bonnet}
Let $F \subset M$ be an admissible surface in a manifold with an
angled polyhedral decomposition. Let $\len(\bdy F \setminus \bdy M)$
be the number of arcs of intersection between $\bdy F \setminus \bdy
M$ and the polyhedra. Then
$$a(F) = -2\pi \chi(F) + 2 \pi \len(\bdy F \setminus \bdy M) \, .$$
\end{prop}

Combinatorial area in angled polyhedra has powerful consequences.
Among them is the following stronger version of a result of Lackenby
\cite[Corollary 4.6]{lack-surg}:

\begin{theorem}\label{hyperbolic-cor}
Let $M$ be an orientable $3$-manifold with an angled polyhedral
decomposition. Then $\bdy M$ is composed of tori, and $M$ is
hyperbolic.
\end{theorem}

\begin{proof}
To prove the first assertion, observe that each component of $\bdy M$
is tiled by boundary faces of the polyhedra. Just inside each boundary
face, a polyhedron has a normal disk of area $0$. These vertex links
glue up to form a closed, boundary-parallel normal surface $F$ of area
$0$.  By Proposition \ref{gauss-bonnet}, $\chi(F) = 0$, and since $M$
is orientable, $F$ must be a torus.

By Thurston's Hyperbolization Theorem \cite{thur-survey}, a manifold
with boundary consisting of tori is hyperbolic if and only if it
contains no essential spheres, disks, tori, or annuli. In our
situation, any such essential surface can be isotoped into normal
form. An essential sphere or disk has positive Euler characteristic,
hence negative area.  Thus it cannot occur.

A normal torus $T \subset M$ has area $0$ and thus, by Lemma
\ref{pos-area}, must be composed of normal disks of area $0$. Since
$T$ has no boundary, these must all be vertex links, which glue up to
form a boundary-parallel torus. Similarly, a normal annulus $A \subset
M$ must be composed entirely of bigons, since a bigon cannot be glued
to a vertex link. But a chain of bigons forms a tube around an edge of
$M$, which is certainly not essential. Thus we can conclude that $M$
is hyperbolic.
\end{proof}

\subsection{Combinatorial length and surgery results}

Lackenby's crucial insight \cite{lack-surg} is that one can use the
combinatorial area of surfaces in a manifold $M$ to define a
combinatorial length of curves on $\bdy M$, and that this notion of
length turns out to be closely related to geometric length on a
maximal cusp.

\begin{define}\label{rel-length-def}
Let $P$ be an angled polyhedron, and let $D \subset P$ be an
admissible disk that intersects at least one boundary face. Let
$\gamma$ be a segment of $\bdy D$ in a boundary face of $P$. Then we
define the {\it length of $\gamma$ relative to $D$} to be
$$\ell(\gamma, D) = \frac{a(D)}{\abs{\bdy D \cap \bdy M}} \, .$$
\end{define}

\begin{define}\label{simplicial-def}
For a manifold $M$ with an angled polyhedral decomposition, let
$\gamma$ be a (possibly non-closed) immersed arc in $\bdy M$. We call 
$\gamma$ a {\it simplicial arc} if
\begin{enumerate}
\item $\gamma$ is disjoint from the vertices of $\bdy M$,
\item the endpoints of $\gamma$ (if any) lie on edges of $\bdy M$,
\item each segment of $\gamma$ in a boundary face is embedded, and
\item no segment of $\gamma$ in a boundary face has endpoints on the
same edge.
\end{enumerate}
\end{define}

We can now define the combinatorial length of simplicial arcs on $\bdy
M$ by considering all the possible {\it inward extensions} of the arc.

\begin{define}\label{length-def}
Let $\gamma \subset \bdy M$ be a simplicial arc.  Let $\gamma_1,
\ldots, \gamma_n$ be the boundary segments that make up $\gamma$,
ordered along a parametrization of $\gamma$.  For each $i$, let $D_i$
be an admissible disk in the corresponding polyhedron, whose boundary
contains $\gamma_i$.  Then $H = \cup_{i=1}^n D_i$ is called an {\it
inward extension} of $\gamma$ if
\begin{enumerate}
\item $\bdy D_i$ agrees with $\bdy D_{i+1}$ on the shared face of their 
polyhedra, and
\item if $\gamma$ is closed, $\bdy D_n$ agrees with $\bdy D_1$ on the 
common face.
\end{enumerate}

We define the {\it combinatorial length} of $\gamma$ to be
$$\lcomb(\gamma) = \inf{\sum_{i=1}^n \ell(\gamma_i, D_i)} \, ,$$
where the infimum is taken over all inward extensions of $\gamma$.
The subscript in $\lcomb(\gamma)$ serves to distinguish combinatorial
length from the geometric cusp length $\lgeom(\gamma)$ used in Section
\ref{geom-cusp}.  When the meaning is clear, we will simply use
$\ell(\gamma)$.
\end{define}

\begin{define}\label{slope-length-def}
Let $s$ be a slope on a boundary component of $M$. Then define
$$\lcomb(s) = \inf{\lcomb(\gamma)} \, ,$$
the infimum being taken over all closed simplicial curves 
$\gamma \subset \bdy M$ that represent non-zero multiples of slope $s$.
\end{define}

The point of this string of definitions is to imply the following
lemma, which is essentially Proposition 4.8 of \cite{lack-surg},
rewritten in terms of polyhedra instead of spines.

\begin{lemma}\label{length-area}
Let $M$ be a manifold with an angled polyhedral decomposition, and let
$F \subset M$ be an admissible surface. Let $C_1, \ldots, C_m$ be the
components of $\bdy F \cap \bdy M$, each $C_j$ representing a non-zero
multiple of some slope $s_{i(j)}$. Then
$$a(F) \geq \sum_{j=1}^{m} \lcomb(s_{i(j)}) \, .$$
\end{lemma}

\begin{proof}
The admissible disks of $F$ bordering on each $C_j$ form one inward
extension of $C_j$. Definition \ref{rel-length-def} has us divide the
area of each disk by the number of its intersections with $\bdy M$, so
we do not end up double-counting any area.
\end{proof}

As a consequence of Lemma \ref{length-area}, surfaces with long
boundary have large combinatorial area, hence large genus. This yields
the following combinatorial analogue of the $6$-Theorem, stated as
Theorem 4.9 of \cite{lack-surg}.

\begin{theorem}[Lackenby]\label{hyplike-surg}
Let $M$ be a manifold with an angled polyhedral decomposition. Let
$s_1, \ldots, s_n$ be a collection of slopes on $\bdy M$, with one
$s_i$ on each component of $\bdy M$. If $\lcomb(s_i) > 2\pi$ for each
$i$, then
the manifold obtained by Dehn filling $M$ along the slopes $s_1,
\ldots, s_n$ is hyperbolike.
\end{theorem}

In fact, Lackenby's machinery allows for an extension of his theorem
to surgeries along only some components of $\bdy M$.

\begin{theorem}\label{partial-surg}
Let $M$ be a manifold with an angled polyhedral decomposition. Let
$s_1, \ldots, s_m$ be a collection of slopes on some, but not all, of the
boundary tori. If $\lcomb(s_i) > 2\pi$ for each $i$, then the
manifold obtained by Dehn filling $M$ along the slopes $s_1,
\ldots, s_m$ is hyperbolic.
\end{theorem}

\begin{proof}
By Thurston's Hyperbolization Theorem \cite{thur-survey}, proving that
the Dehn filled manifold is hyperbolic amounts to ruling out essential
spheres, disks, tori, and annuli.  Any such surface $F$ must intersect
at least one of the solid tori added during the surgery process,
because $M$ is hyperbolic by Theorem \ref{hyperbolic-cor}.  Thus $F$
contains a punctured surface $G \subset M$, whose punctures (not
counting the original boundary components of $F$) represent surgery
slopes $s_{i(1)}, \ldots, s_{i(k)}$ of length greater than $2\pi$.  We
can place $G$ in normal form in the angled polyhedra and compute its
combinatorial area.  Then 
{\setlength\arraycolsep{3pt}
\begin{eqnarray*}
a(G) 
&=& -2\pi\chi(G) \mbox{, \hspace{20pt} by Proposition 
\ref{gauss-bonnet}} \\
&\leq& 2\pi \abs{\bdy G \setminus \bdy F} 
\mbox{, \hspace{9pt} given the choices of $F$} \\
&<& \sum_{j=1}^{\abs{\bdy G \setminus \bdy F}} \ell(s_{i(j)}) 
\mbox{, by assumption} \\
&\leq& a(G) \mbox{, \hspace{40pt} by Lemma \ref{length-area},}
\end{eqnarray*}}
obtaining a contradiction.
\end{proof}

Juxtaposing Theorem \ref{6-theorem} with Theorems \ref{hyplike-surg} and 
\ref{partial-surg}, one can see that for the purpose of ruling out 
exceptional surgeries, $\lcomb(s_i)$ corresponds to $\frac{\pi}{3}
\lgeom(s_i)$.  It turns out that on the cusps of $E(J)$, geometric and
combinatorial length have a similar correspondence (compare Theorem
\ref{thm:curve-lengths} with Corollary \ref{gen-estimate}).  This
yields a second, combinatorial, proof of Theorem \ref{main}.

\section{Normal Surfaces in the Augmented Link Polyhedra}\label{our-polyhedra}

In this section, we apply the normal surface theory of Section
\ref{normal-surf} to the ideal polyhedral decomposition of the
augmented link complement $E(J)$, constructed in Section
\ref{aug-links}.  Recall that by Theorem \ref{newlink-poly-straight},
$P_1$ and $P_2$ are convex ideal polyhedra in $\HH^3$, so they satisfy
the definition of an angled polyhedron.  (See \cite[Theorem
1]{rivin-topology}.)  In fact, they are examples of a special type of
ideal polyhedron, which we call {\it rectangular-cusped}. If we
truncate the ideal vertices, as we did with $P_1$ and $P_2$ in Section
\ref{aug-links}, the resulting {\it boundary faces} subdivide $\bdy M$
into rectangles.

\subsection{Rectangular-cusped polyhedra}

\begin{define}
Let $P$ be an angled ideal polyhedron (see Definition \ref{angled-def})
in which we have truncated the ideal vertices.  We say that $P$ is
{\it rectangular-cusped} if
\begin{enumerate}
\item each boundary face of $P$ (each face of $P \cap \bdy M$) meets $4$ 
interior edges, and
\item each interior edge is labeled with angle $\pi/2$.
\end{enumerate}
\end{define}

Rectangular-cusped polyhedra have two convenient features. First,
their interior faces can be two-colored, in a similar fashion to the
white and shaded faces of $P_1$ and $P_2$.  Around each rectangular
boundary face, opposite interior faces have the same color.  Second,
making all dihedral angles equal to $\pi/2$ ensures that all
combinatorial areas are multiples of $\pi/2$.

In addition to the vertex links and boundary bigons of area $0$ (see
Figure \ref{zero-area-fig}), we need to define a third kind of special
admissible disk.

\begin{figure}[hbt]
\begin{center}
\psfrag{(a)}{(a)}
\psfrag{(b)}{(b)}
\includegraphics{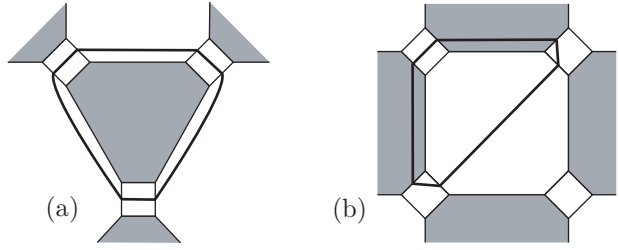}
\caption{Two ideal triangles in a rectangular-cusped polyhedron.}
\label{ideal-triang-face}
\end{center}
\end{figure}

\begin{define}
Let $P$ be a truncated ideal polyhedron. An admissible disk $D \subset
P$ is called an {\it ideal triangle} if 
\begin{enumerate}
\item $\bdy D \subset \bdy P$,
\item $\bdy D$ intersects the boundary faces of $P$ exactly three
times, and 
\item $\bdy D$ is disjoint from the interior edges of $P$. 
\end{enumerate}
\end{define}

Two examples are shown in Figure \ref{ideal-triang-face}. Note that an
ideal triangle $D$ has area $a(D) = \pi$ and length $\ell(\gamma, D) =
\pi/3$ for each segment $\gamma$ of $\bdy D \cap \bdy M$.

\begin{prop}\label{rel-length-bound}
Let $D \subset P$ be an admissible disk in a rectangular-cusped
polyhedron, such that $\bdy D$ passes through at least one boundary
face.  Let $\gamma \subset \bdy M$ be a boundary segment of $\bdy D$.
If $D$ is not a bigon or an ideal triangle,
$$\ell(\gamma, D) \geq \frac{\pi}{2} \, .$$
\end{prop}

\begin{proof}
We consider different cases, conditioned on $n = \abs{\bdy D \cap \bdy
M}$. By Definition \ref{rel-length-def}, $\ell(\gamma, D) = a(D)/n$.

\smallskip

\noindent \underline{\it Case 1: $n=1$.} For this case, we need to 
prove that $a(D) \geq \pi/2$. An admissible disk with one component of
$\bdy D \cap \bdy M$ cannot be a vertex link or boundary bigon, so by
Lemma \ref{pos-area}, $a(D) > 0$. Since all areas in a
rectangular-cusped polyhedron are multiples of $\pi/2$, $a(D) \geq
\pi/2$.

\smallskip

\noindent \underline{\it Case 2: $n=2$.} For this case, we need to 
prove that $a(D) \geq \pi$. If $a(D) = 0$, $D$ is a boundary bigon,
excluded by the hypotheses. So we need to rule out the possibility
that $a(D) = \pi/2$.

If such a disk were to occur, it would have to have $\bdy D \subset
\bdy P$, and $\bdy D$ would have to intersect exactly one interior
edge. Then $\bdy D$ passes through three interior faces, which cannot
all have the same color because two of them share an edge. Thus a
segment $\gamma_1 \subset \bdy D$ in a boundary face must connect
adjacent interior faces, for otherwise all three interior faces would
have the same color. See Figure \ref{case2-fig} for a schematic
picture.

\begin{figure}[b]

\psfrag{g1}{$\gamma_1$}
\psfrag{g2}{$\gamma_2$}
\psfrag{g2'}{$\gamma_2'$}
\psfrag{g2"}{$\gamma_2''$}
\psfrag{g3}{$\gamma_3$}
\psfrag{g4}{$\gamma_4$}
\psfrag{g5}{$\gamma_5$}
\psfrag{g5'}{$\gamma_5'$}
\psfrag{e1}{$e_1$}
\psfrag{e2}{$e_2$}
\psfrag{e3}{$e_3$}
\psfrag{e4}{$e_4$}

\begin{center}
\includegraphics{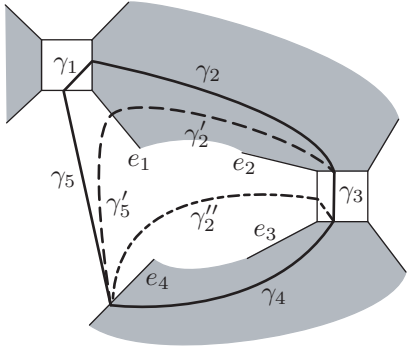}
\caption{Schematic picture for Case $2$ of Proposition
\ref{rel-length-bound}.}
\label{case2-fig}
\end{center}
\end{figure}

We can pull $\gamma_1$ off the boundary face and have it intersect
interior edge $e_1$. This creates a new disk $D'$ with one segment on
$\bdy M$ and area $0$, since this isotopy reduced the area by
$\pi/2$. If $D'$ were admissible, it would be a counterexample to Case
$1$. Thus $\bdy D'$ must violate some condition of admissibility.  The
only way this can happen is if one of the new segments of $\bdy D'$,
$\gamma_2'$ or $\gamma_5'$, has both endpoints on the same edge, or on
adjacent interior and boundary edges.  But since $D$ is admissible,
$e_1$ and $e_4$ must be distinct edges, so $\gamma_5'$ has endpoints
on distinct edges.

Thus $\gamma_2'$ connects adjacent interior and boundary edges, and so
$e_1=e_2$. We can then isotope $\gamma_2'$ across this interior edge,
creating a new disk $D''$ that has just one intersection with $\bdy M$
and one intersection with an interior edge.  Since $\gamma_4$ and the
new segment $\gamma_2''$ lie in adjacent faces of $P$, we have
$e_3=e_4$.  Then $\gamma_4$ connects adjacent interior and boundary
edges, contradicting the assumption that $D$ was admissible.
Therefore, such a disk $D$ does not exist.

\smallskip

\noindent \underline{\it Case 3: $n=3$.} For this case, we need to prove
that $a(D) \geq 3\pi/2$.  The three components of $\bdy D \cap \bdy M$
already ensure that $a(D) \geq \pi$.  So if $\bdy D$ also intersects
an interior edge or the interior of $P$, we have $a(D) \geq 3\pi/2$.
Otherwise, $D$ is an ideal triangle, excluded by the hypotheses.

\smallskip

\noindent \underline{\it Case 4: $n \geq 4$.} For this case,
$$ a(D) \geq n \cdot \pi - 2\pi \geq \frac{n}{2} \cdot \pi \, ,$$
proving the lemma. 
\end{proof}

Thus ideal triangles are the only admissible disks of nonzero area
that contribute less than $\pi/2$ to combinatorial length.  To obtain
the best possible bounds on the length of surgery curves, we need to
find out more about how these triangles fit into polyhedra $P_1$ and
$P_2$ that decompose the link complement $E(J)$.

\subsection{More on ideal triangles} 

\begin{lemma}\label{triang-normal}
Let $P$ be a truncated ideal polyhedron, and let $D \subset P$ be an
ideal triangle. Then all the segments of $\bdy D$ lie in distinct faces of
$\bdy P$, and $D$ is normal.
\end{lemma}

\begin{proof}
$\bdy D$ consists of six segments, alternating between boundary and
interior faces.  Label them $\gamma_1, \ldots, \gamma_6$.  If two of
these segments (say, $\gamma_1$ and $\gamma_3$) lie in the same face
of $P$, then a third segment ($\gamma_2$) must have both endpoints on
the same edge, violating the definition of an admissible disk.  Thus
each $\gamma_i$ lies in a different face, so $\bdy D$ is embedded.
Since $\bdy D$ intersects each boundary face at most once and is
disjoint from the interior edges altogether, $D$ must be normal.
\end{proof}

For the rest of this section, we will work directly with the polyhedra
$P_1$ and $P_2$, and the only manifolds we will consider are $E(J)$
and its Dehn fillings.

\begin{define}\label{triang-types}
In polyhedra $P_1$ and $P_2$, we will classify ideal triangles into
three types. A triangle of {\it type S} is one that is parallel to a
shaded face, as in Figure \ref{ideal-triang-face}(a). A triangle of
{\it type W} is one that is parallel to a white face, as in Figure
\ref{ideal-triang-face}(b).  An ideal triangle parallel to no face of
its polyhedron will be of {\it type N}.
\end{define}

\begin{lemma}\label{triang-classify}
Let $D$ be an ideal triangle in $P \in \{P_1, P_2\}$.  Let $\gamma_1,
\ldots, \gamma_6$ be the segments of $\bdy D$.  Then the following hold:
\begin{enumerate}
\item If $D$ is of type S or type W, then at least two of the $\gamma_i$ 
are parallel to interior edges of $P$.
\item If $D$ is of type N, then no $\gamma_i$ is parallel to an interior 
edge, and the three interior faces of $P$ intersecting $\bdy D$ are
all white faces.
\end{enumerate}
\end{lemma}

\begin{proof}
We consider two cases:

\smallskip
\noindent \underline{\it Case 1: $\bdy D$ intersects both white and 
shaded faces.}  Then $D$ can be schematically represented by the left
side of Figure \ref{triang-proof-fig}.  Label the stumps of interior
edges $e_1, \ldots, e_4$, as in the figure; some of these are likely
to be part of the same edge.  Now, we can pull segment $\gamma_2$ off
the boundary face and have $\bdy D$ intersect edge $e_2$ instead.
This creates a disk $D'$ of area $\pi/2$, which could a priori be
normal.  However, by Case $2$ of Proposition \ref{rel-length-bound},
there are no normal disks that have two intersections with $\bdy M$
and area $\pi/2$.  Thus $D'$ fails some part of Definition
\ref{normal-def}.

\begin{figure}

\psfrag{g1}{$\gamma_1$}
\psfrag{g1'}{$\gamma_1'$}
\psfrag{g2}{$\gamma_2$}
\psfrag{g3}{$\gamma_3$}
\psfrag{g3'}{$\gamma_3'$}
\psfrag{g4}{$\gamma_4$}
\psfrag{g5}{$\gamma_5$}
\psfrag{g6}{$\gamma_6$}
\psfrag{e1}{$e_1$}
\psfrag{e2}{$e_2$}
\psfrag{e3}{$e_3$}
\psfrag{e4}{$e_4$}
\psfrag{e23}{$e_2 = e_3$}
\psfrag{e14}{$e_1 = e_4$}
\psfrag{F}{$F$}
\psfrag{arrow}{$\displaystyle{\Longrightarrow}$}

\begin{center}
\includegraphics{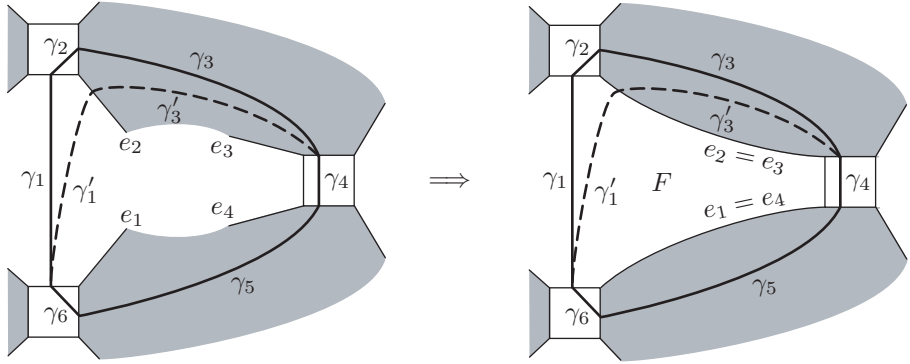}
\caption{Schematic picture of an ideal triangle intersecting faces of both
colors. The shading is generic, and might be reversed.}
\label{triang-proof-fig}
\end{center}
\end{figure}

Since $D$ is normal by Lemma \ref{triang-normal}, the only way that
$D'$ can fail to be normal is if one of the new segments, $\gamma_1'$
or $\gamma_3'$, connects adjacent boundary and interior faces.  If
$\gamma_1'$ violates normality, $e_1$ is the same edge as $e_2$.  But
then $\gamma_3$ and $\gamma_5$ must lie in the same face,
contradicting Lemma \ref{triang-normal}.

If $\gamma_3'$ violates normality by connecting adjacent boundary and
interior faces, we can tighten $\bdy D'$ by removing its intersection
with $e_2=e_3$.  This creates a new disk $D''$ with area $0$.  Segment
$\gamma_5$ and the isotopic image of $\gamma_3$ lie in distinct faces
because they are on opposite sides of edge $e_4$.  So $D''$ is normal,
and thus a boundary bigon.  Then we can conclude that $e_1=e_4$, and
the original disk $D$ was parallel to face $F$, into which we have
pulled $\gamma_3$ (see Figure \ref{triang-proof-fig}, right).  So $D$
is of type S or W.  Notice that both $\gamma_3$ and $\gamma_5$ are
parallel to edges of $F$.

\smallskip
\noindent \underline{\it Case 2: All interior faces intersecting 
$\bdy D$ are the same color.}  If some segment $\gamma_i$ is parallel
to an interior edge, we can isotope $\bdy D$ across that edge, into a
face of a different color, putting us in Case 1.  Otherwise, if no
$\gamma_i$ is parallel to an interior edge, the three interior faces
must all be white.  (Shaded faces are all triangles, in which any arc
connecting distinct ideal vertices is parallel to an edge.)  By Lemma
\ref{triang-normal}, the segments $\gamma_i$ all lie in distinct
faces, so since none of them is parallel to an edge, $D$ cannot be
parallel to a face.  Thus $D$ is of type N, and satisfies conclusion
($2$) of the Lemma.
\end{proof}

\begin{corollary}\label{glue-types}
In an admissible surface in $E(J)$, an ideal triangle of type N
cannot be glued to a bigon or a triangle of type S.
\end{corollary}

\begin{proof}
Let $F$ be a shaded face of $P_1$ or $P_2$, and $D$ be a type S ideal
triangle parallel to $F$.  Since shaded faces are all triangles, every
interior segment of $\bdy D$ is parallel to an interior edge of $F$,
hence an edge of $E(J)$.  Similarly, both interior segments on the
boundary of a bigon are parallel to an edge of $E(J)$.  On the other
hand, by Lemma \ref{triang-classify} the boundary of a type N ideal
triangle does not have any segments parallel to interior edges. 
\end{proof}

\subsection{Progressive arcs and length estimates}

We are now ready to estimate the combinatorial length of surgery
slopes on $\bdy E(J)$.  

\begin{define}
Let $T$ be a torus of $\bdy E(J)$.  Recall that, by Definition
\ref{lattice-basis}, its universal cover $\tild{T}$ contains a lattice 
of shaded and white faces, generated by a basis $\langle \s, \w
\rangle$.  If $T$ is a crossing circle cusp, we will say that the $\w$
direction is {\it meridional} and the $\s$ direction is {\it
longitudinal}. If $T$ is a knot strand cusp, we will say that the $\s$
direction is {\it meridional} and the $\w$ direction is {\it
longitudinal}. (By Lemma \ref{new-cusp-comb}, the meridian and
longitude of $T$ are in fact aligned primarily in these directions.)
\end{define}

Thus if a segment $\gamma$ spans opposite edges of a boundary face $B
\subset \bdy E(J)$, it makes sense to talk of $\gamma$ lying in a
meridional or longitudinal direction.

\begin{define}\label{segment-types}
Let $P \in \{P_1, P_2\}$, and let $D \subset P$ be an admissible disk.
Then $D$ can intersect a boundary face $B \subset \bdy E(J)$ in one of
three types of segments: a {\it longitudinal segment}, connecting
opposite edges of $B$ in a longitudinal direction; a {\it meridional
segment}, connecting opposite edges of $B$ in a meridional direction;
or a {\it diagonal segment}, connecting adjacent edges of $B$.
\end{define}

To estimate the combinatorial length of surgery slopes on $\bdy E(J)$
representing a surgery slope, it helps to divide a curve into smaller
pieces.

\begin{define}\label{prog-arc-def}
Let $T$ be a torus of $\bdy E(J)$, and let $\gamma \subset T$ be a
non-closed simplicial arc (see Definition \ref{simplicial-def}).  Lift
$\gamma$ to an arc $\tild{\gamma} \subset \tild{T}$, and cut
$\tild{T}$ into vertical strips along meridional faces in the
lattice. We say that $\gamma$ is a {\it progressive arc} if
$\tild{\gamma}$ is contained entirely in one of these vertical strips,
and the endpoints of $\tild{\gamma}$ lie on opposite sides of the
strip.
\end{define}

\begin{figure}[h] 
\begin{center}
\includegraphics{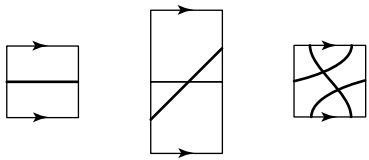}
\caption{The three types of progressive arcs.}
\label{prog-arc-fig}
\end{center}
\end{figure}

In other words, a progressive arc on a crossing circle cusp has
endpoints on consecutive white faces, and constitutes a step in the $\s$
direction. A progressive arc on a knot cusp has endpoints on consecutive
shaded faces, and constitutes a step in the $\w$ direction. In either
case, a progressive arc $\gamma$ can consist of (a) a single
longitudinal segment, (b) two diagonal segments connecting to
different meridians, or (c) two diagonals with some number of
meridional segments between them. These basic types are shown in
Figure \ref{prog-arc-fig}.

\begin{lemma}\label{prog-arc-estimate}
Let $\gamma \subset \bdy E(J)$ be a progressive arc. Then $\ell(\gamma)
\geq \pi/3$.
\end{lemma}

\begin{proof}
Let $H$ be an inward extension of $\gamma$ (see Definition
\ref{length-def}).  For each admissible disk $D_i \subset H$ bordering on
a segment $\gamma_i \subset \gamma$, $\ell(\gamma_i, D_i) = 0$ if and
only if $D_i$ is a boundary bigon.  By Proposition
\ref{rel-length-bound}, every other type of disk contributes at least
$\pi/3$ to $\ell(\gamma)$.  So the only way to have $\ell(\gamma) <
\pi/3$ is if $H$ consists only of bigons.  However, a string of bigons
circles around a single edge of $E(J)$, which means that its
intersection with a component of $\bdy E(J)$ cannot be a progressive
arc.
\end{proof}

\begin{corollary}\label{gen-estimate}
Let $T$ be a torus of $\bdy E(J)$, and let $s$ be a non-trivial
surgery slope on $T$.  If $T$ comes from a crossing circle $C_i$, let
$n$ be the number of crossings in region $R_i$; if $T$ comes from a
component $K_j$ of $K$, let $n$ be the number of twist regions visited
by $K_j$, counted with multiplicity.  Then, in either case,
$$\lcomb(s) \geq \frac{n\pi}{3} \, .$$
\end{corollary}

\begin{proof}
By Theorem \ref{surg-curve-in-basis}, a surgery curve on a crossing
circle corresponding to $n$ crossings must cross at least $n$ (white)
meridional faces, and any surgery curve on a component of $K$ passing
through $n$ twist regions with multiplicity must cross at least $n$
(shaded) meridional faces.  Specifically, we can say that they must
each contain at least $n$ progressive arcs.  Thus the result follows
from Lemma \ref{prog-arc-estimate}.
\end{proof}

For surgery curves on a crossing circle cusp, which by Theorem
\ref{surg-curve-in-basis} look like $n\s \pm \w$ in the basis 
$\langle \s, \w \rangle$, we can obtain a slightly better estimate.

\begin{prop}\label{crossing-circ-estimate}
Let $s \subset \bdy E(J)$ be a surgery slope on a crossing circle cusp
that yields $n$ crossings. Then we have the strict inequality 
$$\lcomb(s) > \frac{n\pi}{3} \, .$$
\end{prop}

\begin{proof}
By Corollary \ref{gen-estimate}, we must only rule out equality.
Equality occurs when a simplicial curve $c$ representing $s$ contains
exactly $n$ progressive arcs, an inward extension of $c$ picks up
length exactly $\pi/3$ per progressive arc, and any part of $c$ not
covered by progressive arcs contributes zero length.  Consider such a
curve.

If a progressive arc $\gamma \subset c$ has combinatorial length
$\pi/3$, it must have an inward extension whose area comes from a
single triangle $D$.  $D$ cannot be of type W, because white faces are
meridional on a crossing circle cusp, and thus a triangle of this
type, plus some bigons, cannot have their boundary segments add up to
a progressive arc.  Thus $D$ must be a triangle of type S or type N.

Let $H$ be an inward extension of $c$.  We claim that if $H$ contains
a type-N triangle, then it consists entirely of type-N triangles.
This is because by Corollary \ref{glue-types}, a type-N triangle $D$
cannot be glued to a type-S triangle or a bigon, and any other type of
admissible disk glued to $D$ would contribute extra area and bring the
total length above $n\pi/3$.  But if $H$ consists entirely of type-N
triangles, $c$ consists entirely of longitudinal segments and never
travels in the $\w$ direction.  Thus we can conclude that $H$ cannot
contain any type-N triangles.

The only remaining possibility is that $H$ consists entirely of type-S
triangles and bigons.  But in this case, all of $H$ is parallel to a
single shaded disk, and again $c$ never traverses the lattice in the
$\w$ direction.  Thus the assumption that $\ell(c) = n\pi/3$ leads to
a contradiction.
\end{proof}

We are now in a position to prove the theorems listed
in the introduction.

\medskip

\noindent {\bf Theorem \ref{hyp-link}.}
{\it Let $K \subset S^3$ be a link with a prime, twist-reduced diagram
$D(K)$. If $D(K)$ has at least two twist regions and every twist
region of $D(K)$ contains at least $6$ crossings, then $K$ is
hyperbolic.  }

\begin{proof}
The assumption that $D(K)$ has at least two twist regions ensures that
the constructions and results of Section \ref{aug-links} apply.  Thus,
by Theorem \ref{newlink-poly-straight}, $K$ is obtained by Dehn
surgery on the crossing circles of a hyperbolic link $J$.  By
Proposition \ref{crossing-circ-estimate}, every surgery slope $s_i$ on
a crossing circle $C_i$ has combinatorial length $\ell(s_i) >
2\pi$. Therefore, by Theorem \ref{partial-surg}, $E(K)$ is hyperbolic.
\end{proof}

\noindent {\bf Theorem \ref{main}.}
{\it Let $K$ be a link in $S^3$ with a prime, twist-reduced diagram
$D(K)$.  Suppose every twist region of $D(K)$ contains at least
$6$ crossings and each component of $K$ passes through at least $7$
twist regions (counted with multiplicity). Then
  \begin{enumerate}
  \item any non-trivial Dehn filling of some but not all components 
  of $K$ is hyperbolic, and
  \item any non-trivial Dehn filling of all the components of $K$ is 
  hyperbolike.
  \end{enumerate}
}

\begin{proof}
By Corollary \ref{gen-estimate}, any non-trivial slope $s$ on a
component of $K$ will have $\ell(s) > 2\pi$, and by Proposition
\ref{crossing-circ-estimate}, the same is true for surgery slopes on
the crossing circles. Thus all surgery slopes on $\bdy E(J)$ are
sufficiently long. Conclusion ($1$) now follows by Theorem
\ref{partial-surg}, and conclusion ($2$) by Theorem
\ref{hyplike-surg}.
\end{proof}

\noindent{\bf Theorem \ref{genus-bound}.}
{\it
  Let $K \subset S^3$ be a link of $k$ components with a prime,
  twist-reduced diagram   $D(K)$. If $D(K)$ has $t \geq 2$ twist
  regions and at least $6$ crossings in each twist region, then
  $$\genus(K) \geq \left\lceil 1 + \frac{t}{6} - \frac{k}{2}
  \right\rceil ,$$
  where $\lceil \cdot \rceil$ is the ceiling
  function that rounds up to the nearest integer.
}

\begin{proof}
  Let $F$ be a Seifert surface for $K$, that is, an orientable incompressible
  surface whose boundary is $K$. Then $F$ contains a punctured surface
  $G \subset E(J)$, where $\bdy G$ consists of curves $\gamma_1,
  \ldots, \gamma_k$ that are longitudes of $K$ and curves
  $\gamma_{k+1}, \ldots \gamma_{k+n}$ along the crossing circles. We
  can place $G$ in normal form in the polyhedra $P_1$ and $P_2$ and
  compute its combinatorial area. Observe that, by Corollary
  \ref{gen-estimate}, the total length of $\gamma_1, \ldots, \gamma_k$
  is at least $2t\pi/3$, because $K$ passes through each twist region
  twice. By Proposition \ref{crossing-circ-estimate}, $\ell(\gamma_i)
  > 2\pi$ for $i > k$. Thus we can compute that
  {\setlength\arraycolsep{3pt}
\begin{eqnarray*}
2\pi \cdot \genus(F) 
&=& 2\pi \cdot \genus(G)  \\
&=& 2\pi \left( 1 - \half \, \chi(G)  - \half (k+n) \right)     \\
&=& 2\pi + \half \, a(G) - \pi k - \pi n \\
&\geq& 2\pi + \half \, \sum_{i=1}^{k} \ell(\gamma_i) - \pi k 
+ \half \sum_{i=k+1}^{k+n} \ell(\gamma_i) - \pi n \\
&\geq& 2\pi + \frac{t\pi}{3} - \pi k \\
&=& 2\pi \left( 1 + \frac{t}{6} - \frac{k}{2} \right).
\end{eqnarray*}}
Since the genus of $F$ is an integer, we are done.
\end{proof}

Observe that the inequality in the computation is an equality whenever
$G$ doesn't meet any crossing circles and consists of only ideal
triangles. This can happen when the twist regions of $D(K)$ always
meet in threes and $G$ lies in the projection plane. In this
situation, Theorem \ref{genus-bound} actually gives the exact value
for the genus of $K$.

\bibliographystyle{hamsplain}

\bibliography{biblio.bib}
\end{document}